\documentclass[11pt,onecoulme]{article}
\usepackage{amsfonts}
\usepackage{mathrsfs}
\usepackage{amssymb}
\usepackage{color, amsmath,amssymb, amsfonts, amstext,amsthm, latexsym}

\usepackage{amssymb, epsfig, amssymb, latexsym}
\usepackage{amsmath}
\usepackage{graphicx}
\usepackage{longtable}
\usepackage{appendix}
\DeclareMathOperator{\esssup}{esssup}
\DeclareMathOperator{\essinf}{essinf}
\numberwithin{equation}{section}

\allowdisplaybreaks

\newtheorem{definition}{Definition}[section]
\newtheorem{theorem}[definition]{Theorem}
\newtheorem{lemma}[definition]{Lemma}
\newtheorem{remark}[definition]{Remark}

\newtheorem{hyp}[definition]{Hypothesis}

  \renewcommand\appendix{\par
    \setcounter{section}{0}
    \setcounter{subsection}{0}
    \gdef\thesection{ Appendix \Alph{section}}}
\title{Viscosity Solutions of Second Order Path-Dependent Partial Differential Equations and Applications 
}
\author{ Shanjian Tang\footnote{Institute of Mathematical Finance and Department of Finance and Control Sciences,
	School of Mathematical Sciences, Fudan University,
	Shanghai 200433, P. R. China, sjtang@fudan.edu.cn. This author is partially supported by National Science Foundation of China (Grant
	No. 11631004) and National Key R\&D Program of China (Grant No. 2018YFA0703903).} \quad Jianjun Zhou\footnote{College of Science, Northwest A$\&$F University, Yangling 712100, Shaanxi, P. R. China, zhoujianjun@nwsuaf.edu.cn.}}

      \date{}

\begin{document}

\maketitle

\pagestyle{plain}

\begin{abstract}
In this article, a notion of viscosity solutions is introduced for fully nonlinear second order path-dependent partial differential equations in the spirit of [Zhou, \emph{Ann. Appl. Probab.,} 33 (2023), 5564-5612].
We prove the existence, comparison principle, consistency and stability for the viscosity solutions. Application to path-dependent stochastic differential games is given.

\medskip

\end{abstract}


\section{Introduction}

 \par
This paper studies viscosity solutions of the following fully nonlinear
path-dependent partial differential equation (PPDE):
 \begin{eqnarray}\label{hjb1}
\begin{cases}
\mathcal{L}V(\gamma_t):=\partial_tV(\gamma_t)+{\mathbf{F}}(\gamma_t,V(\gamma_t),\partial_xV(\gamma_t),\partial_{xx}V(\gamma_t))= 0,\ \ \  (t,\gamma_t)\in
                               [0,T)\times {\Lambda},\\
 V(\gamma_T)=\phi(\gamma_T), \ \ \ \gamma_T\in {\Lambda}_T.
 \end{cases}
\end{eqnarray}
Here,  $\Lambda_t$ is the set of all
 continuous $\mathbb{R}^{d}$-valued functions $\gamma$ defined over $[0,t]$, and  let ${\Lambda}^s=\bigcup_{l\in[s,T]}{\Lambda}_{l}$ and ${\Lambda}$ denote ${\Lambda}^0$; the pathwise (or functional or Dupire; see \cite{dupire1, cotn0, cotn1}) derivatives $\partial_t,\partial_x$ and $\partial_{xx}$ are defined through a functional It\^o formula initiated by
  \cite{dupire1} (see also \cite{cotn0, cotn1}). Such equations arise naturally in many
applications. For example, the dynamic programming equation associated with a
stochastic control problem of non-Markov diffusions (see \cite{ekren3,zhou5}) and the one associated
with a stochastic differential game with non-Markov dynamics (see \cite{pham,zhang}) both fall in
the class of (\ref{hjb1}).

\par
   Studies of   the path-dependant Bellman equation are referred to Peng \cite{peng2},  and  the series of papers:  Ekren, Keller, Touzi
and Zhang \cite{ekren1} and  Ekren,
Touzi and Zhang  
\cite{ekren3, ekren4}, which developed  a  notion of viscosity solutions for the path-dependant Bellman equation in terms of a nonlinear expectation. 
The second author \cite{zhou5} proposed
a notion of viscosity solutions for  the path-dependant Bellman equation in the Crandall-Lions framework, which is further used in this paper, and established its wellposedness: existence, uniqueness, consistency  and stability.
The main innovation of our  approach
is that, due to the lack of regularity  of the supremum  norm $||\cdot||_0$, we  construct a smooth gauge-type functional which is equivalent to $||\cdot||_0^6$.

  Fleming and Souganidis \cite{fle2} first studied  the zero-sum stochastic differential games (SDGs) and  showed that the
lower and upper game values are the viscosity solutions of the corresponding
Bellman-Isaacs equations and coincide under the Isaacs condition. Since then, a lot of
work followed the approach in \cite{fle2}. Buckdahn and Li \cite{buck1} generalized the result in \cite{fle2} with the
help of the theory of backward stochastic differential equations (BSDEs). Pham and Zhang \cite{pham} and  Possamai,  Touzi and   Zhang \cite{pos} studied the
game value of zero-sum path-dependent SDGs (PSDGs) in  weak formulation by the dynamic programming principle and the viscosity solution theory of path-dependent Isaacs equation.
Zhang \cite{zhang} considered  the  zero-sum PSDGs and showed that the existence of   game
values under Isaacs condition  utilizing  a series of approximate state dependent games and their viscosity solution theory.

Our first objective in this paper is to extend  the theory  of Crandall-Lions viscosity solutions to  general fully nonlinear PPDE (\ref{hjb1}).
We adopt the  notion of viscosity solutions introduced in Zhou \cite{zhou5} and give  assumptions to ensure the wellposedness.   Similar to Zhou \cite{zhou5}, we overcome the  difficulty that the supremum norm is not differentiable, and prove the   comparison principle of viscosity solutions. The main difference from Zhou \cite{zhou5} is  that   we prove the existence of viscosity solutions with the Peron's method, while Zhou \cite{zhou5} proved the existence of viscosity solutions with the value function of the optimal control problem.

An alternative  objective of the paper is to apply our results to PSDGs and show that the upper and the
lower value functionals are the unique viscosity solutions of the upper and the lower
path-dependent Hamilton-Jacobi-Bellman-Isaacs equations (PHJBIEs) when the coefficients  are uniformly Lipschitz in the path function under $||\cdot||_0$.
We remark that, Pham and Zhang \cite{pham} proved the uniqueness  only when the diffusion coefficient is uniformly non-degenerate and
the dimension $d$ is either 1 or 2; Possamai,  Touzi and  Zhang \cite{pos}
 did not consider the uniqueness of viscosity solutions; and Zhang \cite{zhang} did not study the existence and uniqueness of  viscosity solutions.
Then none of these results in the above papers are directly applicable to our case.

              The rest of the paper  is organized as follows. In the following
              section, we provide the notations of pathwise derivatives introduced in  \cite{cotn1} and \cite{dupire1}, a modification of Ekeland-Borwein-Preiss variational principle and the smooth gauge-type functions  $\overline{\Upsilon}_0^m$ which are useful in what follows.
           In Section 3, 
            we define classical and viscosity solutions to our
             PPDE (\ref{hjb1}) and  prove  the existence of  viscosity solutions.
             The consistency with the notion of classical solutions and the stability result are also given.
               Section 4  is devoted to proof of   the comparison principle of viscosity solutions to  PPDE (\ref{hjb1}).  In Section 5, we apply our results to PSDGs and show that the upper and the
lower value functionals are the unique viscosity solutions of the upper and the lower PHJBIEs.
\section{Preliminaries}  \label{RDS}
\vbox{}
2.1. \emph{Pathwise derivatives}.
For the vectors $x,y\in \mathbb{R}^{d}$, the scalar product is denoted by $(x,y)_{\mathbb{R}^{d}}$ and the
         Euclidean norm $(x,x)^{\frac{1}{2}}_{\mathbb{R}^{d}}$ is denoted by $|x|$ (we use the same symbol $|\cdot|$ to denote the Euclidean norm on ${\mathbb{{R}}}^k$, for any $k\in {\mathbf{N}}^+$). If $A$ is a vector or a matrix, its transpose is denoted by $A^\top$; For  a matrix $A$, denote its operator norm and
          Hilbert-Schmidt norm    by $|A|$ and $|A|_2$, respectively. Denote by $\mathcal{S}(\mathbb{R}^{d})$  the set of all $(d\times d)$ symmetric matrices.
Let 
$T>0$ be a fixed number.
       For each  $t\in[0,T]$, let $\Lambda_t:= C([0,t],\mathbb{R}^d)$ be the set of all continuous $\mathbb{R}^d$-valued functions defined over $[0,t]$. We denote ${\Lambda}^t=\bigcup_{s\in[t,T]}{\Lambda}_{s}$  and let  ${\Lambda}$ denote ${\Lambda}^0$.  We define a norm on  ${\Lambda}_t$  and a metric on ${\Lambda}$ as follows:
  for any $0\leq t\leq \bar{t}\leq T$ and $\gamma_t,\bar{\gamma}_{\bar{t}}\in {\Lambda}$,
\begin{eqnarray}\label{2.1}
   ||\gamma_t||_{0}:=\sup_{0\leq s\leq t}|\gamma_t(s)|,\ \ d_{\infty}(\gamma_t,\bar{\gamma}_{\bar{t}})= d_{\infty}(\bar{\gamma}_{\bar{t}},\gamma_t):=|t-\bar{t}|
               +||\gamma_{t,\bar{t}}-\bar{\gamma}_{\bar{t}}||_{0},
\end{eqnarray}
where
\begin{eqnarray*}
  \gamma_{t,\bar{t}}(s)&:=&\gamma_t(s){\mathbf{1}}_{[0,t)}(s)+\gamma_t(t){\mathbf{1}}_{[t,\bar{t}]}(s), \ \ s\in [0,\bar{t}].
\end{eqnarray*}
In the sequel,
 for notational simplicity,
we use $||\gamma_{t}-\bar{\gamma}_{\bar{t}}||_{0}$ to denote $||\gamma_{t,\bar{t}}-\bar{\gamma}_{\bar{t}}||_{0}$.
 Then 
 $({\Lambda}_t, ||\cdot||_{0})$ is a Banach space and $({\Lambda}^t, d_{\infty})$ is a complete metric space.
 \par
 \begin{definition}\label{definitionc0409}
       Let $t\in[0,T)$ and $f:{\Lambda}^t\rightarrow K$ be given for some Hilbert space $K$.
\begin{description}
        \item{(i)}
                 We say $f\in C^0({\Lambda}^t,K)$ (resp., $f\in USC^0({\Lambda}^t,K)$, $f\in LSC^0({\Lambda}^t,K)$) if $f$ is continuous (resp., upper semicontinuous, lower semicontinuous) in $\gamma_s$  on ${\Lambda}^t$ under $d_{\infty}$.
                  \item{(ii)}
                 We say $f\in C_p^0({\Lambda}^t,K)\subset C^0({\Lambda}^t,K)$ if $f$ grows  in a polynomial way.
\end{description}
\end{definition}
For notational simplicity, we abbreviate  $C^0({\Lambda}^t,\mathbb{R})$, $USC^0({\Lambda}^t,\mathbb{R})$, $LSC^0({\Lambda}^t,\mathbb{R})$ and $C_p^0({\Lambda}^t,\mathbb{R})$ as $C^0({\Lambda}^t)$,  $USC^0({\Lambda}^t)$,  $LSC^0({\Lambda}^t)$ and $C^0_p({\Lambda}^t)$, respectively.
 \par
 We now define the path derivatives via the functional It\^o formula,
which is initiated by Dupire \cite{dupire1}, and plays a crucial role in this  paper.
 \begin{definition}\label{theoremito}
Let $t\in[0,T)$  be given.  We say $f\in C^{1,2}_p(\Lambda^t)$ if $f\in C^0_p(\Lambda^t)$ and there exist $\partial_tf\in C^0_p(\Lambda^t)$, $\partial_xf\in C^0_p(\Lambda^t,\mathbb{R}^d)$, $\partial_{xx}f\in C^0_p(\Lambda^t,{\cal{S}}(\mathbb{R}^d))$
such that, for any  continuous adapted process $X$ on $[0,T]$ which is a  semi-martingale  on $[t,T]$,
\begin{eqnarray}\label{statesop}
                f(X_s)=f(X_t)+\int_{t}^{s}\partial_tf(X_l)dl+\frac{1}{2}\int^{s}_{t}\partial_{xx}f(X_l)d\langle X\rangle(l)
                +\int^{s}_{t}\partial_xf(X_l)dX(l),\   s\in[t,T], \ \mathbb{P}\mbox{-}a.s.\
\end{eqnarray}
 Here and in the following, for every $s\in [0,T]$, $X(s)$ denotes  the value of $X$  at
 time $s$, and $X_s$ the whole history path of $X$ from time 0 to $s$.
\end{definition}
\par
We remark that the above $\partial_tf, \partial_xf, \partial_{xx}f$, if they exist, are unique. In fact,
for  every $(\hat{t},\gamma_{\hat{t}})\in [t,T)\times \Lambda^{\hat{t}}$, $\alpha\in \mathbb{R}^{d}$ and $\beta\in \mathbb{R}^{d\times n}$,
 let  $$
           X(s)=\gamma_{\hat{t}}({\hat{t}})+\int^{s}_{{\hat{t}}}\alpha dl+\int^{s}_{{\hat{t}}}\beta dW(l),\ \ s\in [{\hat{t}},T],
 $$
 and $X(s)=\gamma_{\hat{t}}(s)$, $s\in [0,{\hat{t}})$, where $\{W(t),t \geq 0\}$ be a $n$-dimensional standard Wiener process. Then $X(\cdot)$ is a continuous semi-martingale on $[t,T]$. First, let $\alpha=\mathbf{0}$, $\beta=\mathbf{0}$, together with the required regularity
  $\partial_tf\in C^0_p(\Lambda^t)$ we get
  \begin{eqnarray}\label{2.3}
               \partial_tf(\gamma_{\hat{t}}):=\lim_{h\rightarrow0,h>0}\frac{1}{h}\left[f(\gamma_{{\hat{t}},{\hat{t}}+h})-f(\gamma_{\hat{t}})\right],\ ({\hat{t}},\gamma_{\hat{t}})\in [t,T)\times{\Lambda}.
\end{eqnarray}
Here and
in the sequel, for notational simplicity, we use $\mathbf{0}$ to denote the element or the function  which is
identically equal to zero.
Next, let $\beta=\mathbf{0}$,   by the continuity of $\partial_x{f}$  and the arbitrariness of $\alpha\in {\mathbb{R}^d}$, we  have the uniqueness of $\partial_xf$. Finally,  let   $\alpha=\mathbf{0}$,
by the continuity and symmetry  of $\partial_{xx}{f}$ and the arbitrariness of $\beta\in \mathbb{R}^{d\times n}$, we see that $\partial_{xx}{f}$ is also unique.
\par

\vbox{}
2.2. \emph{ Ekeland-Borwein-Preiss variational principle}. In this subsection, we introduce a modification of Ekeland-Borwein-Preiss variational principle (see Theorem 2.5.2 in  Borwein \& Zhu  \cite{bor1}) which
 plays a crucial role in the proof of the  comparison principle.
We
firstly recall the definition of a gauge-type function for the specific set $\Lambda^t$.
\begin{definition}\label{gaupe}
              Let $t\in [0,T]$ be fixed.  We say that a continuous functional $\lambda:\Lambda^t\times \Lambda^t\rightarrow [0,+\infty)$ is a {gauge-type function} on $\Lambda^t$ provided that:
             \begin{description}
        \item{(i)} $\lambda(\gamma_s,\gamma_s)=0$ for all $(s,\gamma_s)\in [t,T]\times \Lambda^t$,
        \item{(ii)} we have $$\lim_{\lambda(\gamma_s,\eta_l)\to 0} d_\infty(\gamma_s,\eta_l)=0. $$
        \end{description}
\end{definition}
\begin{lemma}\label{theoremleft}  (see  Lemma 2.13 in Zhou \cite{zhou5})    
Let $t\in [0,T]$ be fixed and let $f:\Lambda^t\rightarrow \mathbb{R}$ be an upper semicontinuous functional  bounded from above. Suppose that $\lambda$ is a gauge-type function on $\Lambda^t$
 and $\{\delta_i\}_{i\geq0}$ is a sequence of positive number, and suppose that $\varepsilon>0$ and $(t_0,\gamma^0_{t_0})\in [t,T]\times \Lambda^t$ satisfy
 $$
f(\gamma^0_{t_0})\geq \sup_{(s,\gamma_s)\in [t,T]\times \Lambda^t}f(\gamma_s)-\varepsilon.
 $$
 Then there exist $(\hat{t},\hat{\gamma}_{\hat{t}})\in [t,T]\times \Lambda^t$ and a sequence $\{(t_i,\gamma^i_{t_i})\}_{i\geq1}\subset [t,T]\times \Lambda^t$ such that
  \begin{description}
        \item{(i)} $\lambda(\gamma^0_{t_0},\hat{\gamma}_{\hat{t}})\leq \frac{\varepsilon}{\delta_0}$,  $\lambda(\gamma^i_{t_i},\hat{\gamma}_{\hat{t}})\leq \frac{\varepsilon}{2^i\delta_0}$ and $t_i\uparrow \hat{t}$ as $i\rightarrow\infty$,
        \item{(ii)}  $f(\hat{\gamma}_{\hat{t}})-\sum_{i=0}^{\infty}\delta_i\lambda(\gamma^i_{t_i},\hat{\gamma}_{\hat{t}})\geq f(\gamma^0_{t_0})$, and
        \item{(iii)}  $f(\gamma_s)-\sum_{i=0}^{\infty}\delta_i\lambda(\gamma^i_{t_i},\gamma_s)
            <f(\hat{\gamma}_{\hat{t}})-\sum_{i=0}^{\infty}\delta_i\lambda(\gamma^i_{t_i},\hat{\gamma}_{\hat{t}})$ for all $(s,\gamma_s)\in [\hat{t},T]\times \Lambda^{\hat{t}}\setminus \{(\hat{t},\hat{\gamma}_{\hat{t}})\}$.

        \end{description}
\end{lemma}

\par

\vbox{}
2.3. \emph{Functionals $\Upsilon^{m}$}.    In this subsection we introduce the  functionals $\Upsilon^{m}$, which are the key to proving  the comparison principle and  stability  of viscosity solutions.
\par
For every $m\in {\mathbf{N}}^+$, introduce 
 \begin{eqnarray*}
&&\Upsilon^m(\gamma_t):=
            \frac{(||\gamma_{t}||_{0}^{2m}-|\gamma_{t}(t)|^{2m})^3}{||\gamma_{t}||^{4m}_{0}}\mathbf{1}_{\{||\gamma_{t}||_{0}\neq0\}}+3|\gamma_{t}(t)|^{2m}, \ (t,\gamma_t)\in [0,T]\times{{\Lambda}};\\
&&
          \Upsilon^{m}_0(\gamma_t,\eta_s)
          :=\Upsilon^m(\gamma_{t,t\vee s}-\eta_{s,t\vee s}), \ \ \ (t,\gamma_t), (s,\eta_s)\in [0,T]\times{\Lambda};\\
&&
        \overline{\Upsilon}^{m}_0(\gamma_t,\eta_s)
        := \Upsilon^{m}_0(\gamma_t,\eta_s)+|s-t|^2, \ \ \ (t,\gamma_t), (s,\eta_s)\in [0,T]\times{\Lambda}.
\end{eqnarray*}

For simplicity, we let  $\Upsilon$, $\Upsilon_0$ and $\overline{\Upsilon}_0$ denote $\Upsilon^3$, $\Upsilon^{3}_0$ and $\overline{\Upsilon}^{3}_0$, respectively.
Combining   Theorem 2.3 and  Lemmas 3.1, 3.2 and 3.4 in \cite{zhou5}, we have
\begin{lemma}\label{theoremS}
For any integer $m\ge 2$,
 $\Upsilon^{m}(\cdot)\in C^{1,2}_p({\Lambda})$  and
   \begin{eqnarray}\label{0528a}
   \partial_t\Upsilon^{m}(\gamma_t)=0;
   \end{eqnarray}
 \begin{eqnarray}\label{0528a}
    \partial_{x}\Upsilon^{m}(\gamma_t)=6m\left(1-\frac{(||\gamma_t||^{2m}_{0}-|\gamma_t(t)|^{2m})^2}{||\gamma_t||^{4m}_{0}}\right)|\gamma_t(t)|^{2m-2}\gamma_t(t)\mathbf{1}_{\{||\gamma_t||_{0}\neq0\}};
\end{eqnarray}
  \begin{eqnarray}\label{0528b}
    \partial_{xx}\Upsilon^{m}(\gamma_t)
    &=&\bigg{[}\frac{24m^2(||\gamma_t||_{0}^{2m}-|\gamma_t(t)|^{2m})|\gamma_t(t)|^{4m-4}
    \gamma_t(t)(\gamma_t(t))^\top}
   {{||\gamma_t||_{0}^{4m}}}\nonumber\\
   &&+12m(m-1)\left(1-\frac{(||\gamma_t||_{0}^{2m}-|\gamma_t(t)|^{2m})^2}{{||\gamma_t||^{4m}}}\right)
   |\gamma_t(t)|^{2m-4}\gamma_t(t)(\gamma_t(t))^\top\nonumber\\
&&+6m\left(1-\frac{(||\gamma_t||_{0}^{2m}-|\gamma_t(t)|^{2m})^2
  }{{||\gamma_t||_{0}^{4m}}}\right)|\gamma_t(t)|^{2m-2}I\bigg{]}\mathbf{1}_{\{||\gamma_t||_{0}\neq0\}}. \ \
\end{eqnarray}
  Moreover, the following estimates hold: for any $(t,\gamma_t,\eta_t)\in [0,T]\times {\Lambda}\times {\Lambda}$,
\begin{eqnarray}\label{s0}
                  |\partial_{x}\Upsilon^{m}(\gamma_t)|\leq 6m|\gamma_t(t)|^{2m-1}, \ \  |\partial_{xx}\Upsilon^{m}(\gamma_t)|\leq 6m(6m-1)|\gamma_t(t)|^{2m-2};
\end{eqnarray}
and
\begin{eqnarray}\label{up}
||\gamma_t||_{0}^{2m}\leq   \Upsilon^{m}(\gamma_t)
                \leq 3||\gamma_t||_{0}^{2m}, \ \ \left(\Upsilon^{m}(\gamma_t+\eta_t)\right)^{\frac{1}{2m}}\leq \left( \Upsilon^{m}(\gamma_t)\right)^{\frac{1}{2m}}+ \left(\Upsilon^{m}(\eta_t)\right)^{\frac{1}{2m}}.
\end{eqnarray}
Finally, for every $t\in [0,T]$, $\overline{\Upsilon}_0(\cdot,\cdot)$ is a gauge-type function on compete metric space $(\Lambda^t,d_{\infty})$.
\end{lemma}

\section{Viscosity solutions to  PPDEs: Existence.}

\par
In this section, we consider the  second order path-dependent partial differential equation (PPDE) (\ref{hjb1}). As usual, we start with classical solutions. For every fixed $(\hat{t},\hat{\gamma}_{\hat{t}})\in [0,T)\times \Lambda$ and $r\in \mathbb{R}^+\cup\{\infty\}$, define
$$B_r(\hat{\gamma}_{\hat{t}}):=\{\eta_s:s\geq \hat{t}, \ \ d_{\infty}(\eta_s,\hat{\gamma}_{\hat{t}})<r\}.$$
It is clear that $\Lambda^{\hat{t}}=B_{\infty}(\hat{\gamma}_{\hat{t}})$ for every $(\hat{t},\hat{\gamma}_{\hat{t}})\in [0,T)\times \Lambda$.
\par
\begin{definition}\label{definitionccc}     (Classical solution) For every fixed $(\hat{t},\hat{\gamma}_{\hat{t}})\in [0,T)\times \Lambda$ and $r\in \mathbb{R}^+\cup\{\infty\}$,
              A functional $v\in C_p^{1,2}({\Lambda}^{\hat{t}})$       is called a classical solution  (resp. subsolution, supersolution) to the PPDE (\ref{hjb1}) on $B_r(\hat{\gamma}_{\hat{t}})$ if the terminal condition, $v(\gamma_T)=(\mbox{resp}., \leq, \geq)\phi(\gamma_T)$ for all $\gamma_T\in B_r(\hat{\gamma}_{\hat{t}})$ is satisfied, and
              $$
                                  \mathcal{L}v(\gamma_t)=(\mbox{resp.}\ \geq,\ \leq)0,\ \ \ \forall\ (t,\gamma_t)\in [0,T)\times B_r(\hat{\gamma}_{\hat{t}}).
              $$
 \end{definition}
In what follows, by a $modulus \ of \ continuity$, we mean a continuous function $\rho:[0,\infty)\rightarrow[0,\infty)$, with $\rho(0)=0$ and  subadditivity: $\rho(t+s)\leq \rho(t)+\rho(s)$, for all $t,s>0$; by a $local\ modulus\ of \  continuity$, we mean  a continuous function $\rho:[0,\infty)\times[0,\infty) \rightarrow[0,\infty)$, with the properties that, for each $r\geq0$, $t\rightarrow \rho(t,r)$ is a modulus of continuity and $\rho$ is non-decreasing in second variable.
\par
We will make the following assumptions about the function $\mathbf{F}:\Lambda\times \mathbb{R}\times \mathbb{R}^d\times {\cal{S}}(\mathbb{R}^d)\rightarrow\mathbb{R}$.
        \begin{hyp}\label{hypstate}
\begin{description}
        \item{(i)}
        $\mathbf{F}$ is 
          continuous on 
            $\Lambda\times \mathbb{R}\times \mathbb{R}^d\times {\cal{S}}(\mathbb{R}^d)$.
\par
       \item{(ii)}
                 There exists a constant $\nu\geq0$
                 such that, for every $(t,\gamma_t,r,p,X)
                 \in [0,T]\times\Lambda\times  \mathbb{R}\times \mathbb{R}^d\times {\cal{S}}(\mathbb{R}^d)$,
      \begin{eqnarray*}
              \mathbf{F}(\gamma_t,r,p,X)-\mathbf{F}(\gamma_t,s,p,X)\geq
                 \nu(s-r) \ \  \mbox{when} \  r\leq s.
\end{eqnarray*}
\item{(iii)}  For every $(t,\gamma_t,r,p)\in [0,T]\times\Lambda\times  \mathbb{R}\times \mathbb{R}^d$
      \begin{eqnarray*}
              \mathbf{F}(\gamma_t,r,p,X)\leq\mathbf{F}(\gamma_t,r,p,Y) \ \  \mbox{when} \  X\leq Y.
\end{eqnarray*}
\item{(iv)} 
There exists a  local modulus of  continuity $\rho$ such that, for every $(t,\gamma_t,\eta_t,r)\in [0,T]\times\Lambda\times\Lambda\times  \mathbb{R}$,  
    for any
$\beta>0$, for all $X,Y\in {\cal{S}}(\mathbb{R}^d)$ satisfying
 \begin{eqnarray*}
                              -3\beta\left(\begin{array}{cc}
                                    I&0\\
                                    0&I
                                    \end{array}\right)\leq \left(\begin{array}{cc}
                                    X&0\\
                                    0&Y
                                    \end{array}\right)\leq  3\beta \left(\begin{array}{cc}
                                    I&-I\\
                                    -I&I
                                    \end{array}\right),
\end{eqnarray*}
 we have
      \begin{eqnarray*}
              &&\mathbf{F}(\gamma_t,r,\beta(\gamma_t(t)-\eta_t(t)),X)-\mathbf{F}(\eta_t,r,\beta(\gamma_t(t)-\eta_t(t)),-Y)\\
              &\leq& \rho(\beta||\gamma_t-\eta_t||^2_0+||\gamma_t-\eta_t||_0,|r|\vee||\gamma_t||_0\vee||\eta_t||_0).
\end{eqnarray*}
\item{(v)} There exists a constant $M_F\geq 0$  such that
      \begin{eqnarray*}
              |\mathbf{F}(\gamma_t,r,p+q,X+Y)-\mathbf{F}(\gamma_t,r,p,X)|
              \leq M_F[(1+||\gamma_t||_0)|q|+(1+||\gamma_t||_0^2)|Y|].
\end{eqnarray*}
\end{description}
\end{hyp}
       Now we turn to viscosity solutions.
                      For every $(t,\gamma_t)\in [0,T]\times \Lambda$ and $w\in C^0(\Lambda)$, define
$$
             \mathcal{A}^+(\gamma_t,w):=\Big{\{}\varphi\in C^{1,2}_p({\Lambda}^{t}):  0=({w}-{{\varphi}})({\gamma_t})=\sup_{(s,\eta_s)\in [t,T]\times{\Lambda}}
                         ({w}- {{\varphi}})(
                         \eta_s)\Big{\}},
$$
and
$$
             \mathcal{A}^-(\gamma_t,w):=\Big{\{}\varphi\in C^{1,2}_p({\Lambda}^{t}):  0=({w}-{{\varphi}})({\gamma_t})=\inf_{(s,\eta_s)\in [t,T]\times{\Lambda}}
                         ({w}-{{\varphi}})(
                         \eta_s)\Big{\}}.
$$
\begin{definition}\label{definition4.1} \ \
 $w\in USC^0({\Lambda})$ (resp., $w\in LSC^0({\Lambda})$) is called a
                             viscosity subsolution (resp.,  supersolution)
                             to  (\ref{hjb1}) if the terminal condition,  $w(\gamma_T)\leq \phi(\gamma_T)$(resp.,  $w(\gamma_T)\geq \phi(\gamma_T)$) for all
                             $\gamma_T\in {\Lambda}_T$ is satisfied,
                                and whenever  ${\varphi}\in \mathcal{A}^+(\gamma_s,w)$ (resp.,  ${\varphi}\in \mathcal{A}^-(\gamma_s,w)$)  with $(s,{\gamma}_{s})\in [0,T)\times{\Lambda}$,  we have
\begin{eqnarray*}
                         \mathcal{L}\varphi(\gamma_s)\geq0  \ \  (\mbox{resp.},\ \mathcal{L}\varphi(\gamma_s)\leq0).
\end{eqnarray*}
                                $w\in C^0({\Lambda})$ is said to be a
                             viscosity solution to PPDE (\ref{hjb1}) if it is
                             both a viscosity subsolution and a viscosity
                             supersolution.
\end{definition}
\par
We are now in a  position  to give  the existence  result 
for the viscosity solutions.
\begin{theorem}\label{theoremvexist12} \ \
                          Let Hypothesis \ref{hypstate} (i)  be satisfied. Let $\cal{A}$ be a family of viscosity subsolution of (\ref{hjb1}). Suppose that
                          there exist a  local modulus of  continuity $\rho$ and  a constant $\Delta>0$  such that, for every $(t,\gamma_t, w)\in [0,T]\times\Lambda\times{\cal{A}}$, 
                        {  \begin{eqnarray}\label{202109250}
                          w(\gamma_t)\leq w(\gamma_{t,s})+\rho(|s-t|,||\gamma_t||_0), \ \  
                           s\in [t,T\wedge (t+\Delta)].
                          \end{eqnarray}}
                          Let
\begin{eqnarray}\label{20210925}
                          u(\gamma_t):=\sup\{w(\gamma_t): w\in {\cal{A}}\},\ \ (t,\gamma_t)\in [0,T]\times\Lambda
\end{eqnarray}
and assume that $u^*(\gamma_t)<\infty$ for all $(t,\gamma_t)\in [0,T]\times\Lambda$.
 Then $u^*$ is a viscosity subsolution of (\ref{hjb1}).
   Here $u^*$ is the upper
semicontinuous envelope of $u$ (see    
 \cite[Definition D.10]{fab1}), i.e.,
$$
u^*(\gamma_t)=\limsup_{(s,\eta_s)\in [0,T]\times \Lambda, (s,\eta_s)\rightarrow(t,\gamma_t)}u(\eta_s).
$$
Similarly, $u_*$ is
the lower semicontinuous envelope of $u$, i.e.,
$$
u_*(\gamma_t)=\liminf_{(s,\eta_s)\in [0,T]\times \Lambda, (s,\eta_s)\rightarrow(t,\gamma_t)}u(\eta_s).
$$
\end{theorem}
{\bf  Proof}. \ \
First, for every ${\gamma}_{T}\in \Lambda_T$, by (\ref{202109250}), there exists a sequence $(\gamma^n_{T}, u_n)\in \Lambda_T\times {\cal{A}}$ such that
              \begin{eqnarray*}
                     (\gamma^n_{T}, u_n(\gamma^n_{T}))\rightarrow ({\gamma}_{T},u^*({\gamma}_{T}))\ \mbox{as}\ n\rightarrow\infty.
                    \end{eqnarray*}
Since $u_n\in {\cal{A}}$, we have
$$
       u_n(\gamma^n_{T})\leq \phi(\gamma^n_{T}).
$$
Letting $n\rightarrow\infty$,
$$
       u^*({\gamma}_{T})\leq \phi({\gamma}_{T}).
$$
 Next, let  $\varphi\in \mathcal{A}^+(\hat{\gamma}_{\hat{t}},u^*)$
                  with
                   $(\hat{t},\hat{\gamma}_{\hat{t}})\in [0,T)\times \Lambda$.
                       By (\ref{202109250}), there is a sequence $(t_n,\gamma^n_{t_n}, u_n)\in [\hat{t},T)\times \Lambda\times {\cal{A}}$ such that
              \begin{eqnarray}\label{0424a}
                     (t_n,\gamma^n_{t_n}, u_n(\gamma^n_{t_n}))\rightarrow (\hat{t},\hat{\gamma}_{\hat{t}},u^*(\hat{\gamma}_{\hat{t}})) \ \mbox{as}\ n\rightarrow\infty.
                    \end{eqnarray}
                     Set
 $$\Gamma_n(\gamma_t):=(u_n-\varphi)(\gamma_t)-\overline{\Upsilon}_0(\gamma_t,\hat{\gamma}_{\hat{t}}),\ \ (t,\gamma_t)\in [\hat{t},T]\times \Lambda^{\hat{t}}.$$
 Then, the functional $\Gamma_n$ is upper semicontinuous.      By $\varphi\in \mathcal{A}^+(\hat{\gamma}_{\hat{t}},u^*)$ and the definition of $u^*$,
         \begin{eqnarray}\label{0424b}
         \Gamma_n(\gamma_t)\leq (u_n-\varphi)(\gamma_t)\leq (u^*-\varphi)(\gamma_t)\leq 0,\ \   (t,\gamma_t)\in [\hat{t},T]\times \Lambda^{\hat{t}}.
         \end{eqnarray}
        This means that $\Gamma_n$ is bounded from above on $\Lambda^{\hat{t}}$.
               Define a sequence of positive numbers $\{\delta_i\}_{i\geq0}$  by 
        $\delta_i=\frac{1}{2^i}$ for all $i\geq0$.  Since $\overline{\Upsilon}_0(\cdot,\cdot)$ is a gauge-type function on compete metric space $(\Lambda^{\hat{t}},d_{\infty})$, for every $n>0$ and $\varepsilon>0$,
           from Lemma \ref{theoremleft} it follows that,
  for every  $(\check{t}_{0},\check{\gamma}^{0}_{\check{t}_{0}})\in [\hat{t},T]\times  \Lambda^{\hat{t}}$ satisfying
\begin{eqnarray}\label{20210925a}
\Gamma_n(\check{\gamma}^{0}_{\check{t}_{0}})
\geq \sup_{(t,\gamma_t)\in [\hat{t},T]\times \Lambda^{\hat{t}}}\Gamma_n(\gamma_t)-\frac{1}{n}\geq \Gamma_n(\gamma^n_{t_n})-\frac{1}{n},\
\end{eqnarray}
  there exist $(\hat{t}_n,\hat{\gamma}^{n}_{\hat{t}_{n}})\in [\hat{t},T]\times \Lambda^{\hat{t}}$ and sequence $\{(\check{t}_{i},\check{\gamma}^{i}_{\check{t}_{i}})\}_{i\geq1} \subset
  [\hat{t},T]\times \Lambda^{\hat{t}}$ such that
  \begin{description}
        \item{(i)} $\overline{\Upsilon}_0(\check{\gamma}^{0}_{\check{t}_{0}},\hat{\gamma}^{n}_{\hat{t}_{n}})\leq \frac{1}{n}$,
         $\overline{\Upsilon}_0(\check{\gamma}^{i}_{\check{t}_{i}},\hat{\gamma}^{n}_{\hat{t}_{n}})
        \leq \frac{1}{2^in}$
          and $\check{t}_{i}\uparrow \hat{t}_{n}$ as $i\rightarrow\infty$,
        \item{(ii)}  $\Gamma_n(\hat{\gamma}^{n}_{\hat{t}_{n}})
            -\sum_{i=0}^{\infty}\frac{1}{2^i}\overline{\Upsilon}_0(\check{\gamma}^{i}_{\check{t}_{i}},\hat{\gamma}^{n}_{\hat{t}_{n}})
        \geq \Gamma_{n}(\check{\gamma}^{0}_{\check{t}_{0}})$, and
        \item{(iii)}    for all $(t,\gamma_t)\in [\hat{t}_{n},T]\times \Lambda^{\hat{t}_{n}}\setminus \{(\hat{t}_{n},\hat{\gamma}^{n}_{\hat{t}_{n}})\}$,
        \begin{eqnarray*}
        \Gamma_{n}(\gamma_t)
        -\sum_{i=0}^{\infty}
        \frac{1}{2^i}\overline{\Upsilon}_0(\check{\gamma}^{i}_{\check{t}_{i}},\gamma_t)
            <\Gamma_{n}(\hat{\gamma}^{n}_{\hat{t}_{n}})
            -\sum_{i=0}^{\infty}\frac{1}{2^i}\overline{\Upsilon}_0(\check{\gamma}^{i}_{\check{t}_{i}},\hat{\gamma}^{n}_{\hat{t}_{n}}).
        \end{eqnarray*}
        \end{description}
By (\ref{0424b}), we have $u_n-\varphi\leq 0$ on  $\Lambda^{\hat{t}}$.  Then by (\ref{20210925a}) and the property (ii) of $(\hat{t}_n,\hat{\gamma}^{n}_{\hat{t}_{n}})$,
\begin{eqnarray}\label{0424c}
&&-\overline{\Upsilon}_0(\hat{\gamma}^{n}_{\hat{t}_{n}},\hat{\gamma}_{\hat{t}})\geq(u_n-\varphi)(\hat{\gamma}^{n}_{\hat{t}_{n}})-\overline{\Upsilon}_0(\hat{\gamma}^{n}_{\hat{t}_{n}},\hat{\gamma}_{\hat{t}})
=\Gamma_n(\hat{\gamma}^{n}_{\hat{t}_{n}})\nonumber\\
&\geq& \Gamma_n(\gamma^n_{t_n})-\frac{1}{n}=(u_n-\varphi)(\gamma^n_{t_n})-\overline{\Upsilon}_0(\gamma^n_{t_n},\hat{\gamma}_{\hat{t}})-\frac{1}{n}.
 \end{eqnarray}
 Notice that, by (\ref{0424a}), $\varphi\in \mathcal{A}^+(\hat{\gamma}_{\hat{t}},u^*)$, (\ref{up}) and the definition of $\overline{\Upsilon}_0$,
 $$
               (u_n-\varphi)(\gamma^n_{t_n})\rightarrow (u^*-\varphi)(\hat{\gamma}_{\hat{t}})=0\ \mbox{and}\ \overline{\Upsilon}_0(\gamma^n_{t_n},\hat{\gamma}_{\hat{t}})\rightarrow0 \ \mbox{as}\ n\rightarrow\infty.
 $$
Letting $n\rightarrow\infty$ in (\ref{0424c}), 
  we get
\begin{eqnarray}\label{04072}
         \lim_{n\to \infty}  \overline{\Upsilon}_0(\hat{\gamma}^{n}_{\hat{t}_{n}},\hat{\gamma}_{\hat{t}})=0.
\end{eqnarray}
Since $\hat{t}<T$, then for sufficiently large integers $n$,
 $$\varphi_1:=\varphi+\overline{\Upsilon}_0(\cdot,\hat{\gamma}_{\hat{t}})+\sum_{i=0}^{\infty}
        \frac{1}{2^i}\overline{\Upsilon}_0(\check{\gamma}^{i}_{\check{t}_{i}},\cdot)\quad \in \quad \mathcal{A}^+(\hat{\gamma}^{n}_{\hat{t}_{n}},u_n)$$
                  with
                   $(\hat{t}_n,\hat{\gamma}^{n}_{\hat{t}_{n}})\in [0,T)\times \Lambda$.
  Since $u_n$ is a viscosity subsolution of (\ref{hjb1}), we have
\begin{eqnarray}\label{04071}
                          \partial_t{\varphi_{1}}(\hat{\gamma}^{n}_{\hat{t}_{n}})
                           +{\mathbf{F}}(\hat{\gamma}^{n}_{\hat{t}_{n}}, u_n(\hat{\gamma}^{n}_{\hat{t}_{n}}),
                           \partial_x{\varphi_{1}}(\hat{\gamma}^{n}_{\hat{t}_{n}}),\partial_{xx}{\varphi_{1}}(\hat{\gamma}^{n}_{\hat{t}_{n}}))\geq0.
\end{eqnarray}
Notice that
$$
\partial_t{\varphi_{1}}(\gamma_t)=\partial_t{\varphi}(\gamma_t)+2(t-\hat{t})+2\sum_{i=0}^{\infty}
        \frac{1}{2^i}(t-\check{t}_{i}),
$$
$$
\partial_x{\varphi_{1}}(\gamma_t)=\partial_x{\varphi}(\gamma_t)+\partial_x\Upsilon(\gamma_t-\hat{\gamma}_{\hat{t},t})+\partial_x\left[\sum_{i=0}^{\infty}
        \frac{1}{2^i}\Upsilon(\gamma_t-\check{\gamma}^{i}_{\check{t}_{i},t})\right],
$$$$
\partial_{xx}{\varphi_{1}}(\gamma_t)=\partial_{xx}{\varphi}(\gamma_t)+\partial_{xx}\Upsilon(\gamma_t-\hat{\gamma}_{\hat{t},t})+\partial_{xx}\left[\sum_{i=0}^{\infty}
        \frac{1}{2^i}\Upsilon(\gamma_t-\check{\gamma}^{i}_{\check{t}_{i},t})\right],
$$
and, by (\ref{s0}) and  the property (i) of $(\check{t}_{i},\check{\gamma}^{i}_{\check{t}_{i}})$,
$$
                     |\hat{t}_{n}-\hat{t}|\leq\overline{\Upsilon}^{\frac{1}{2}}_0(\hat{\gamma}^{n}_{\hat{t}_{n}},\hat{\gamma}_{\hat{t}}),
$$
$$
|\partial_x\Upsilon(\hat{\gamma}^{n}_{\hat{t}_{n}}-\hat{\gamma}_{\hat{t},\hat{t}_{n}})|\leq 18|\hat{\gamma}^{n}_{\hat{t}_{n}}(\hat{t}_{n})-\hat{\gamma}_{\hat{t}}(\hat{t})|^5\leq 18\overline{\Upsilon}^{\frac{5}{6}}_0(\hat{\gamma}^{n}_{\hat{t}_{n}},\hat{\gamma}_{\hat{t}}),
$$
$$
|\partial_{xx}\Upsilon(\hat{\gamma}^{n}_{\hat{t}_{n}}-\hat{\gamma}_{\hat{t},\hat{t}_{n}})|\leq 306|\hat{\gamma}^{n}_{\hat{t}_{n}}(\hat{t}_{n})-\hat{\gamma}_{\hat{t}}(\hat{t})|^4\leq 306\overline{\Upsilon}_0^{\frac{2}{3}}(\hat{\gamma}^{n}_{\hat{t}_{n}},\hat{\gamma}_{\hat{t}}),
$$
$$
\left|\sum_{i=0}^{\infty}
        \frac{1}{2^i}(\hat{t}_{n}-\check{t}_{i})\right|\leq \sum_{i=0}^{\infty}
        \frac{1}{2^i}\left(\frac{1}{2^in}\right)^{\frac{1}{2}}=2\left(\frac{1}{n}\right)^{\frac{1}{2}},
$$
$$
\left|\partial_x\left[\sum_{i=0}^{\infty}
        \frac{1}{2^i}\Upsilon(\hat{\gamma}^{n}_{\hat{t}_{n}}-\check{\gamma}^{i}_{\check{t}_{i},\hat{t}_{n}})\right]\right|\leq 18\sum_{i=0}^{\infty}
        \frac{1}{2^i}|\hat{\gamma}^{n}_{\hat{t}_{n}}(\hat{t}_{n})-\check{\gamma}^{i}_{\check{t}_{i}}(\check{t}_{i})|^5\leq 18\sum_{i=0}^{\infty}
        \frac{1}{2^i}\left(\frac{1}{2^in}\right)^{\frac{5}{6}}=36\left(\frac{1}{n}\right)^{\frac{5}{6}},
$$
$$
\left|\partial_{xx}\left[\sum_{i=0}^{\infty}
        \frac{1}{2^i}\Upsilon(\hat{\gamma}^{n}_{\hat{t}_{n}}-\check{\gamma}^{i}_{\check{t}_{i},\hat{t}_{n}})\right]\right|\leq 306\sum_{i=0}^{\infty}
        \frac{1}{2^i}|\hat{\gamma}^{n}_{\hat{t}_{n}}(\hat{t}_{n})-\check{\gamma}^{i}_{\check{t}_{i}}(\check{t}_{i})|^4\leq306\sum_{i=0}^{\infty}
        \frac{1}{2^i}\left(\frac{1}{2^in}\right)^{\frac{2}{3}}=612\left(\frac{1}{n}\right)^{\frac{2}{3}}.
$$
Letting $n\rightarrow\infty$ in (\ref{04071}), by (\ref{0424a}),  (\ref{04072}) and Hypothesis \ref{hypstate} (i), we have
$$
                          \partial_t{\varphi}(\hat{\gamma}_{\hat{t}})
                           +{\mathbf{F}}(\hat{\gamma}_{\hat{t}}, u^*(\hat{\gamma}_{\hat{t}}),
                           \partial_x{\varphi}(\hat{\gamma}_{\hat{t}}),\partial_{xx}{\varphi}(\hat{\gamma}_{\hat{t}}))\geq0.
$$
Thus we show that
 $u^*$ is  a viscosity subsolution to (\ref{hjb1}).\ \ $\Box$
 \begin{theorem}\label{theoremvexist1} \ \
                          Let Hypothesis \ref{hypstate} (i) and (iii)  be satisfied and  comparison hold for (\ref{hjb1}); i.e., if $w$ is a subsolution of (\ref{hjb1}) and $v$ is a
                            supersolution of (\ref{hjb1}), then $w\leq v$. Suppose $\underline{u}$ and $\overline{u}$ be respectively a viscosity subsolution and a viscosity supersolution of (\ref{hjb1}) such that 
                             $\underline{u}(\gamma_T)=\overline{u}(\gamma_T), \ \gamma_T\in \Lambda_T$ and  $\overline{u}^*(\gamma_t)<\infty$, $\underline{u}_*(\gamma_t)>-\infty$ for all $(t,\gamma_t)\in [0,T]\times\Lambda$, and
                             suppose that there exist a  local modulus of  continuity $\rho$ and a constant $\Delta>0$ such that, for every $(t,\gamma_t)\in [0,T]\times\Lambda$, 
                        {  \begin{eqnarray}\label{202109250123}
                          |\underline{u}(\gamma_t)- \underline{u}(\gamma_{t,s})|\vee |\overline{u}(\gamma_t)- \overline{u}(\gamma_{t,s})|\leq\rho(|s-t|, ||\gamma_t||_0), \ \  s\in [t,T\wedge (t+\Delta)].
                          \end{eqnarray}}
                          Then the function
\begin{eqnarray}\label{20210925}
                          u(\gamma_t)=\sup\{w(\gamma_t):& \underline{u}\leq w\leq \overline{u}, w\ \mbox{is a viscosity subsolution of (\ref{hjb1}) and satisfies (\ref{202109250123})}\nonumber\\
                         &\mbox{for the same  local modulus of  continuity }\ \rho \ \mbox{and  constant}\  \Delta>0\}
\end{eqnarray}
  is a viscosity solution of (\ref{hjb1}).
\end{theorem}
{\bf  Proof}. \ \ It is clear that $\underline{u}_*\leq u_*\leq u\leq u^* \leq \overline{u}^*$.  By Theorem \ref{theoremvexist12}  $u^*$ is a viscosity subsolution of (\ref{hjb1}) and hence, by comparison, $u^*\leq \overline{u}$.
It then follows from the definition of $u$ that $u=u^*$ (so $u$ is a viscosity subsolution).  If the condition for $u_*$ being a viscosity supersolution of (\ref{hjb1}) is violated at $\hat{\gamma}_{\hat{t}}$ with $\hat{t}\in [0,T)$
for the test function $\psi$ then
$$
     \psi\in \mathcal{A}^-(\hat{\gamma}_{\hat{t}},u_*) \ \mbox{and}\ \   \partial_t\psi(\hat{\gamma}_{\hat{t}})+{\mathbf{F}}(\hat{\gamma}_{\hat{t}},\psi(\hat{\gamma}_{\hat{t}}),\partial_x\psi(\hat{\gamma}_{\hat{t}}),\partial_{xx}\psi(\hat{\gamma}_{\hat{t}}))>0,
$$
and by continuity,
$u_{\delta,\alpha}(\eta_s)=\delta+\psi(\eta_s)-\alpha\overline{\Upsilon}_0(\eta_s,\hat{\gamma}_{\hat{t}})$
 is a classical
subsolution of (\ref{hjb1}) on $B_r(\hat{\gamma}_{\hat{t}})$ 
for all small $r,\alpha>0$ and $0<\delta\leq \alpha(T-\hat{t})^2$. Since
$$
               u(\eta_s)\geq u_*(\eta_s)\geq \psi(\eta_s),\ \ (s,\eta_s)\in [\hat{t},T]\times\Lambda,
$$
if we choose $\delta=\left(\frac{r^6}{4^6}\wedge \frac{r^2}{4^2}\wedge (T-\hat{t})^2\right)\alpha$ then $u(\eta_s)>u_{\delta,\alpha}(\eta_s)$ for $s\in[\hat{t},T]$ and $d_{\infty}(\eta_s,\hat{\gamma}_{\hat{t}})>\frac{r}{2}$, 
 and then, by Lemma \ref{theorem3.20407jia}, the function
\begin{eqnarray}\label{0407acc}
U(\eta_s)=\begin{cases}
           \max\{{u(\eta_s), u_{\delta,\alpha}(\eta_s)}\},  \ \mbox{if} \   \eta_s\in B_r(\hat{\gamma}_{\hat{t}});\\ 
u(\eta_s),\ \ \ \ \ \ \ \ \ \ \ \ \ \ \ \ \ \  \ \ \   \mbox{otherwise}
\end{cases}
\end{eqnarray}
 is a viscosity subsolution of (\ref{hjb1}). 
   We observe that 
  there are points such that $U(\eta_s)>u(\eta_s)$; in fact, by definition of $u_*$, there is a sequence $(s_n,\eta^n_{s_n})\in [\hat{t},T]\times \Lambda$ such that $(\eta^n_{s_n},u(\eta^n_{s_n}))\rightarrow(\hat{\gamma}_{\hat{t}}, u_*(\hat{\gamma}_{\hat{t}}))$ as $n\rightarrow\infty$, and then
 $$
          \lim_{n\rightarrow\infty}(U(\eta^n_{s_n})-u(\eta^n_{s_n})){\geq} u_{\delta,\alpha}(\hat{\gamma}_{\hat{t}})-u_*(\hat{\gamma}_{\hat{t}})= u_*(\hat{\gamma}_{\hat{t}})+\delta-u_*(\hat{\gamma}_{\hat{t}})>0.
 $$
Clearly, $\underline{u}\leq U$. By comparison, $U\leq \overline{u}$ and since $u$ is the maximal subsolution between $\underline{u}$ and $\overline{u}$, we arrive at the contradiction $U\leq u$. Therefore, $u_*$ is a
 supersolution of (\ref{hjb1}) and then, by comparison for (\ref{hjb1}), $u^*=u\leq u_*$, showing that $u$ is continuous and is a viscosity solution of (\ref{hjb1}). \ \  $\Box$
\par
To complete the proof of Theorem \ref{theoremvexist1}, it remains to state and prove the following two
lemmas. Before doing that, we give the following definition.
\begin{definition}\label{definition4.10407} \ \ For every fixed $(\hat{t},\hat{\gamma}_{\hat{t}})\in [0,T)\times \Lambda$ and $r>0$,
 $w\in USC^0({\Lambda}^{\hat{t}})$  is called a
                             viscosity subsolution
                             to  PPDE (\ref{hjb1}) on $B_r(\hat{\gamma}_{\hat{t}})$ if the terminal condition,  $w(\gamma_T)\leq \phi(\gamma_T)$ for all
                             $\gamma_T\in B_r(\hat{\gamma}_{\hat{t}})$ is satisfied,
                                and whenever  ${\varphi}\in \mathcal{A}^+(\gamma_s,w)$   with $(s,{\gamma}_{s})\in [\hat{t},T)\times B_r(\hat{\gamma}_{\hat{t}})$,  we have
\begin{eqnarray*}
                      \mathcal{L}\varphi(\gamma_s)\geq0.
\end{eqnarray*}
\end{definition}
\begin{lemma}\label{theorem3.20407}
                     Fix $(\hat{t},\hat{\gamma}_{\hat{t}})\in [0,T)\times \Lambda$ and $r>0$. Let Hypothesis \ref{hypstate}  (iii) hold true and  $v\in C_p^{1,2}({\Lambda}^{\hat{t}})$ be a classical subsolution of PPDE (\ref{hjb1}) on $B_r(\hat{\gamma}_{\hat{t}})$. Then
                 $v$ is a viscosity subsolution  of  PPDE (\ref{hjb1}) on $B_r(\hat{\gamma}_{\hat{t}})$.
\end{lemma}
The proof is rather standard and is postponed to Appendix A.
\begin{lemma}\label{theorem3.20407jia}
                      Let Hypothesis \ref{hypstate} (iii) hold true. Then the function $U$ defined in (\ref{0407acc})
                is a viscosity subsolution  of  PPDE (\ref{hjb1}).
\end{lemma}
{\bf  Proof}. \ \  First, from $u=u^*\in  USC^0({\Lambda})$ and $u_{\delta,\alpha}\in C^{1,2}_p(\Lambda^{\hat{t}})$, it follows that $U\in  USC^0({\Lambda})$.
  Second, since  $u$ is a viscosity subsolution of PPDE (\ref{hjb1}) and  $u_{\delta,\alpha}$ is a classical subsolution of PPDE (\ref{hjb1}) on $B_r(\hat{\gamma}_{\hat{t}})$, we have
  $$
           u(\gamma_T)\leq \phi(\gamma_T), \ \gamma_T\in \Lambda_T, \ \mbox{and}\  u_{\delta,\alpha}(\gamma_T)\leq \phi(\gamma_T), \ \gamma_T\in B_r(\hat{\gamma}_{\hat{t}}).
  $$
  Then, by  the definition of $U$, the terminal condition $U(\gamma_T)\leq \phi(\gamma_T), \ \gamma_T\in \Lambda_T$ is satisfied.
 Third,  let   $\varphi\in \mathcal{A}^+({\gamma}_{{t}},U)$
                  with
                   $({t},{\gamma}_{{t}})\in [0,T)\times \Lambda$.
 If $\gamma_t\notin B_r(\hat{\gamma}_{\hat{t}})$, by the definition of $U$,
 \begin{eqnarray}\label{0407acf}
                         0=(u-\varphi)(\gamma_t)=(U-\varphi)(\gamma_t)=\sup_{(s,\eta_s)\in [t,T]\times\Lambda^t}
                         (U- \varphi)(\eta_s)\geq\sup_{(s,\eta_s)\in [t,T]\times\Lambda^t}
                         (u- \varphi)(\eta_s).
\end{eqnarray}
 Since $u$ is a viscosity subsolution of PPDE (\ref{hjb1}), we have
 \begin{eqnarray}\label{0407acd}
                        \mathcal{L}\varphi(\gamma_s)\geq0.
\end{eqnarray}
  If $\gamma_t\in B_r(\hat{\gamma}_{\hat{t}})$ and $u(\gamma_t)\geq u_{\delta,\alpha}(\gamma_t)$, by the definition of $U$, we also have (\ref{0407acf}) and (\ref{0407acd}).\\
  If $\gamma_t\in B_r(\hat{\gamma}_{\hat{t}})$ and $u(\gamma_t)< u_{\delta,\alpha}(\gamma_t)$, noting that $u(\eta_s)>u_{\delta,\alpha}(\eta_s)$ if $s\in[\hat{t},T]$ and $ d_{\infty}(\eta_s,\hat{\gamma}_{\hat{t}})>\frac{r}{2}$, by the definition of $U$,
 $$
                         0=(u_{\delta,\alpha}-\varphi)(\gamma_t)=(U-\varphi)(\gamma_t)=\sup_{(s,\eta_s)\in [t,T]\times\Lambda^t}
                         (U- \varphi)(\eta_s)\geq\sup_{(s,\eta_s)\in [t,T]\times\Lambda^t}
                         (u_{\delta,\alpha}- \varphi)(\eta_s).
$$
 Since $u_{\delta,\alpha}$ is a classical subsolution of PPDE (\ref{hjb1}) on $B_r(\hat{\gamma}_{\hat{t}})$, by Lemma \ref{theorem3.20407},
                 $u_{\delta,\alpha}$ is a viscosity subsolution  of  PPDE (\ref{hjb1}) on $B_r(\hat{\gamma}_{\hat{t}})$ and (\ref{0407acd}) holds true.
                 \par
From above all, the function $U$ defined in (\ref{0407acc})
                is a viscosity subsolution  of  PPDE (\ref{hjb1}).  \ \ $\Box$

   \par
We conclude this section with   the  consistency and stability of viscosity solutions.

\begin{theorem}\label{theorem3.2}
                      Let Hypothesis \ref{hypstate} (iii) hold true and  $v\in C_p^{1,2}({\Lambda})$. Then
                 $v$ is a classical solution (resp.  subsolution, supersolution) of  PPDE (\ref{hjb1}) on $\Lambda$ if and only if it is a viscosity solution  (resp.  subsolution, supersolution) of  PPDE (\ref{hjb1}).
\end{theorem}
\begin{theorem}\label{theoremstability}
                      Let $\mathbf{F}$ satisfy   Hypothesis \ref{hypstate} (i), and $v\in C^0(\Lambda)$. Assume
                       \item{(i)}      for any $\varepsilon>0$, there exist $\mathbf{F}^\varepsilon$ and $v^\varepsilon\in C^0(\Lambda)$ such that  $\mathbf{F}^\varepsilon$ satisfy  Hypothesis \ref{hypstate} (i) and $v^\varepsilon$ is a viscosity solution (resp. subsolution, supersolution) of PPDE (\ref{hjb1}) with generator $\mathbf{F}^\varepsilon$;
                           \item{(ii)} as $\varepsilon\rightarrow0$, $(\mathbf{F}^\varepsilon,v^\varepsilon)$ converge to
                           $(\mathbf{F}, v)$  uniformly in the following sense: 
\begin{eqnarray}\label{sss}
                          \lim_{\varepsilon\rightarrow0}\sup_{(\gamma_t,x,y,z)\in \Lambda\times \mathbb{R}\times \mathbb{R}^d\times {\cal{S}}(\mathbb{R}^d)}[(\mathbf{F}^\varepsilon-\mathbf{F})(\gamma_t,x,y,z)|+|(v^\varepsilon-v)({\gamma}_t)|]=0.
\end{eqnarray}
                   Then $v$ is a viscosity solution (resp. subsolution, supersolution) of  PPDE (\ref{hjb1}) with generator $\mathbf{F}$.
\end{theorem}
The proof of Theorems \ref{theorem3.2} and  \ref{theoremstability}   is rather standard, and is given in the appendix A for  the convenience of the reader.

\section{Viscosity solutions to  PPDEs: Comparison  principle.}
\par
 \vbox{}
             This section is devoted to proof of  the comparison principle of  viscosity
                   solutions to (\ref{hjb1}), which is 
                    required in Theorem \ref{theoremvexist1}. 
The main result of this section is stated as follows.
\begin{theorem}\label{theoremhjbm}
Suppose Hypothesis \ref{hypstate}   holds.
                         Let $W_1\in C^0({\Lambda})$ $(\mbox{resp}., W_2\in C^0({\Lambda}))$ be  a viscosity subsolution (resp., supersolution) to PPDE (\ref{hjb1}) and  let  there exist  constant $L>0$, $m\geq 3$
                          and a local modulus of continuity $\rho$,
                        such that, for any  
                        $(t,\gamma_t),(s,\eta_s)\in[0,T]\times{\Lambda}$,
\begin{eqnarray}\label{w}
                                   |W_1(\gamma_t)|\vee |W_2(\gamma_t)|\leq L (1+||\gamma_t||^m_0);
                                   \end{eqnarray}
\begin{eqnarray}\label{w1}
                               &&|W_1(\gamma_t)-W_1(\eta_s)|\vee|W_2(\gamma_t)-W_2(\eta_s)|\nonumber\\
                               &\leq&
                        \rho(|s-t|,||\gamma_t||_0\vee||\eta_s||_0)+L (1+||\gamma_t||^m_0+||\eta_s||^m_0)||\gamma_{t}-\eta_s||_0.
\end{eqnarray}
                   Then  $W_1\leq W_2$.
\end{theorem}
\par
  The proof follows from the analysis in \cite{zhou5}.   
 We  note that for $\varrho>0$, the functional
                    defined by $\tilde{W}:=W_1-\frac{\varrho}{t+1}$ is a viscosity subsolution
                   for
 \begin{eqnarray}\label{0619c}
\begin{cases}
{\partial_t} \tilde{W}(\gamma_t)+{\mathbf{F}}(\gamma_t, \tilde{W}(\gamma_t), \partial_x \tilde{W}(\gamma_t),\partial_{xx} \tilde{W}(\gamma_t))
          = \frac{\varrho}{(t+1)^2}, \ \  (t,\gamma_t)\in [0,T)\times {\Lambda}, \\
\tilde{W}(\gamma_T)=\phi(\gamma_T), \ \  \gamma_T\in \Lambda_T.
\end{cases}
\end{eqnarray}
                As $W_1\leq W_2$ follows from $\tilde{W}\leq W_2$ in
                the limit $\varrho\downarrow0$, it suffices to prove
                $W_1\leq W_2$ under the following additional assumption:
$$
{\partial_t} {W_1}(\gamma_t)+{\mathbf{F}}(\gamma_t, {W_1}(\gamma_t), \partial_x {W_1}(\gamma_t),\partial_{xx}{W_1}(\gamma_t))
          \geq c,\ \ c:=\frac{\varrho}{(T+1)^2}, \ \  (t,\gamma_t)\in [0,T)\times {\Lambda}.
$$
\par
   {\bf  Proof of Theorem \ref{theoremhjbm}}. \ \
 We only need to prove that $W_1(\gamma_t)\leq W_2(\gamma_t)$ for all $(t,\gamma_t)\in
[T-\bar{a},T)\times
       {\Lambda}$.
        Here,
        $$\bar{a}=\frac{1}{144m^2M_F}\wedge{T}.$$
         Then, we can  repeat the same procedure for the case
        $[T-i\bar{a},T-(i-1)\bar{a})$.  Thus, we assume the converse result that $(\tilde{t},\tilde{\gamma}_{\tilde{t}})\in (T-\bar{a},T)\times
      {\Lambda}$ exists  such that
        $\tilde{m}:=W_1(\tilde{\gamma}_{\tilde{t}})-W_2(\tilde{\gamma}_{\tilde{t}})>0$. 
\par
         Consider that  $\varepsilon >0$ is  a small number such that
 $$
 W_1(\tilde{\gamma}_{\tilde{t}})-W_2(\tilde{\gamma}_{\tilde{t}})-2\varepsilon \frac{\nu T-\tilde{t}}{\nu
 T}\Upsilon^{m}(\tilde{\gamma}_{\tilde{t}})
 >\frac{\tilde{m}}{2},
 $$
      and
\begin{eqnarray}\label{5.3}
                          \frac{\varepsilon}{\nu T}\leq\frac{c}{4},
\end{eqnarray}
             where
$$
            \nu=1+\frac{1}{144m^2M_FT}.
$$
 Next, let $\Lambda^t\otimes\Lambda^t:=\{(\gamma_s,\eta_s)|\gamma_s,\eta_s\in \Lambda^t\}$ for all $ t\in [0,T]$, we define for any $\beta>0$ and  $(\gamma_t,\eta_t)\in {\Lambda}^{T-\bar{a}}\otimes{\Lambda}^{T-\bar{a}}$,
\begin{eqnarray*}
                 \Psi(\gamma_t,\eta_t)=W_1(\gamma_t)-W_2(\eta_t)-{\beta}\Upsilon_0(\gamma_{t},\eta_{t})-\beta^{\frac{1}{3}}|\gamma_{t}(t)-\eta_{t}(t)|^2
                 -\varepsilon\frac{\nu T-t}{\nu
                 T}(\Upsilon^{m}(\gamma_t)+\Upsilon^{m}(\eta_t)).
\end{eqnarray*}
By  Step 1-Step 3 in the proof of Theorem 4.1 in \cite{zhou5},
  for every  $(\gamma^0_{t_0},\eta^0_{t_0})\in  \Lambda^{\tilde{t}}\otimes \Lambda^{\tilde{t}}$ satisfying
$$
\Psi(\gamma^0_{t_0},\eta^0_{t_0})\geq \sup_{(s,(\gamma_s,\eta_s))\in [\tilde{t},T]\times (\Lambda^{\tilde{t}}\otimes \Lambda^{\tilde{t}})}\Psi(\gamma_s,\eta_s)-\frac{1}{\beta},\
\    \mbox{and} \ \ \Psi(\gamma^0_{t_0},\eta^0_{t_0})\geq \Psi(\tilde{\gamma}_{\tilde{t}},\tilde{\gamma}_{\tilde{t}}) >\frac{\tilde{m}}{2},
 $$
  there exist $(\hat{t},(\hat{\gamma}_{\hat{t}},\hat{\eta}_{\hat{t}}))\in [\tilde{t},T]\times (\Lambda^{\tilde{t}}\otimes \Lambda^{\tilde{t}})$ and a sequence $\{(t_i,(\gamma^i_{t_i},\eta^i_{t_i}))\}_{i\geq1}\subset
  [\tilde{t},T]\times (\Lambda^{\tilde{t}}\otimes \Lambda^{\tilde{t}})$ such that
  \begin{description}
        \item{(i)} $\Upsilon_0(\gamma^0_{t_0},\hat{\gamma}_{\hat{t}})+\Upsilon_0(\eta^0_{t_0},\hat{\eta}_{\hat{t}})+|\hat{t}-t_0|^2\leq \frac{1}{\beta}$,
         $\Upsilon_0(\gamma^i_{t_i},\hat{\gamma}_{\hat{t}})+\Upsilon_0(\eta^i_{t_i},\hat{\eta}_{\hat{t}})+|\hat{t}-t_i|^2
         \leq \frac{1}{\beta2^i}$ and $t_i\uparrow \hat{t}$ as $i\rightarrow\infty$,
        \item{(ii)}  $\Psi_1(\hat{\gamma}_{\hat{t}},\hat{\eta}_{\hat{t}})
\geq \Psi(\gamma^0_{t_0},\eta^0_{t_0})$, and
        \item{(iii)}    for all $(s,(\gamma_s,\eta_s))\in [\hat{t},T]\times (\Lambda^{\hat{t}}\otimes \Lambda^{\hat{t}})\setminus \{(\hat{t},(\hat{\gamma}_{\hat{t}},\hat{\eta}_{\hat{t}}))\}$,
        \begin{eqnarray}\label{iii4}
        \Psi_1(\gamma_s,\eta_s)
           <\Psi_1(\hat{\gamma}_{\hat{t}},\hat{\eta}_{\hat{t}}),
        \end{eqnarray}
        \end{description}
        where
$$
     \Psi_1(\gamma_t,\eta_t):=  \Psi(\gamma_t,\eta_t)
        -\sum_{i=0}^{\infty}
        \frac{1}{2^i}[{\Upsilon}_0(\gamma^i_{t_i},\gamma_t)+{\Upsilon}_0(\eta^i_{t_i},\eta_t)+|t-t_i|^2], \ \  \  (\gamma_t,\eta_t)\in  \Lambda^{\tilde{t}}\otimes \Lambda^{\tilde{t}}.
$$
Moreover, there exist ${{M}_0}, N>0$ independent of $\beta$
    such that
                   \begin{eqnarray}\label{5.10jiajiaaaa}
                   ||\hat{\gamma}_{\hat{t}}||_0\vee||\hat{\eta}_{\hat{t}}||_0<M_0, \ \mbox{and} \ \hat{t}\in [\tilde{t},T) \ \ \mbox{for all}\ \ \beta\geq N.
                    \end{eqnarray}
We also have  the following result  holds true:
 \begin{eqnarray}\label{5.10}
                      \beta ||\hat{{\gamma}}_{{\hat{t}}}-\hat{{\eta}}_{{\hat{t}}}||_{0}^6
                         +\beta|\hat{{\gamma}}_{{\hat{t}}}(\hat{t})-\hat{{\eta}}_{{\hat{t}}}(\hat{t})|^4
                         \rightarrow0 \ \mbox{as} \ \beta\rightarrow\infty.
 \end{eqnarray}
             We should note that the point
             $({\hat{t}},\hat{{\gamma}}_{{\hat{t}}},\hat{{\eta}}_{{\hat{t}}})$ depends on $\beta$ and
              $\varepsilon$.
\par
                  We define,
  for $(t,\gamma_t,\eta_t)\in [0,T]\times {\Lambda}\times {\Lambda}$,
\begin{eqnarray}\label{06091}
                             w_{1}(\gamma_t)=W_1(\gamma_t)-2^5\beta \Upsilon_0(\gamma_{t},\hat{\xi}_{\hat{t}})-\varepsilon\frac{\nu T-t}{\nu
                 T}\Upsilon^{m}(\gamma_t)
                 -\varepsilon \overline{\Upsilon}_0(\gamma_t,\hat{\gamma}_{\hat{t}})
               -\sum_{i=0}^{\infty}
        \frac{1}{2^i}\overline{\Upsilon}_0(\gamma^i_{t_i},\gamma_t),
        \end{eqnarray}
        \begin{eqnarray}\label{06092}
                             w_{2}(\eta_t)=-W_2(\eta_t)-2^5\beta \Upsilon_0(\eta_{t},\hat{\xi}_{\hat{t}})-\varepsilon\frac{\nu T-t}{\nu
                 T}\Upsilon^{m}(\eta_t)
                 -\varepsilon \overline{\Upsilon}_0(\eta_t,\hat{\eta}_{\hat{t}})-\sum_{i=0}^{\infty}
        \frac{1}{2^i}{\Upsilon}_0(\eta^i_{t_i},\eta_t),
\end{eqnarray}
   where $\hat{\xi}_{\hat{t}}=\frac{\hat{\gamma}_{\hat{t}}+\hat{\eta}_{\hat{t}}}{2}$.   We  note that $w_1,w_2$ depend on $\hat{\xi}_{{\hat{t}}}$, and thus on $\beta$ and
              $\varepsilon$.
By  Step 4 in the proof of Theorem 4.1 in \cite{zhou5},
  there exist
sequences  $(l_{k},\check{\gamma}^{k}_{l_k}), (s_{k},\check{\eta}^{k}_{s_k})\in [\hat{t},T]\times \Lambda^{\hat{t}}$ and
 the sequences of functionals $(\varphi_k,\psi_k)\in C_p^{1,2}(\Lambda^{l_k})\times C_p^{1,2}(\Lambda^{s_k})$ bounded from below 
    such that
   \begin{eqnarray}\label{4.23}
\lim_{k\rightarrow\infty}[d_\infty(\check{\gamma}^{k}_{l_k},\hat{\gamma}_{\hat{t}})
+d_\infty(\check{\eta}^{k}_{s_k},\hat{\eta}_{\hat{t}})]=0,
\end{eqnarray}
the functional
\begin{eqnarray}\label{0609a}
 w_{1}(\gamma_t)-\varphi_k(\gamma_t), \quad \gamma_t\in \Lambda^{l_k}
\end{eqnarray}
has a strict global maximum $0$ at  $\check{\gamma}^{k}_{l_k}$, while the functional
\begin{eqnarray}\label{0609b}
w_{2}(\eta_t)-\psi_k(\eta_{t}), \quad \eta_t\in \Lambda^{s_k}
\end{eqnarray}
has a strict global maximum $0$ at $\check{\eta}^{k}_{s_k}$, and
\begin{eqnarray}\label{0608v1}
       \lim_{k\to \infty }\begin{pmatrix}
        \partial_t\varphi_k( \check{\gamma}^{k}_{l_{k}}),&\partial_x\varphi_k( \check{\gamma}^{k}_{l_{k}}),&\partial_{xx}\varphi_k( \check{\gamma}^{k}_{l_{k}})\end{pmatrix}=
      \begin{pmatrix}
     b_1,& 2\beta^{\frac{1}{3}} (\hat{{\gamma}}_{{\hat{t}}}({{\hat{t}}})-\hat{{\eta}}_{{\hat{t}}}({{\hat{t}}})), &X\end{pmatrix},
\end{eqnarray}
\begin{eqnarray}\label{0608vw1}
     \lim_{k\to \infty}  \begin{pmatrix}  
       \partial_t\psi_k( \check{\eta}^{k}_{s_{k}}), &\partial_x\psi_k( \check{\eta}^{k}_{s_{k}}), &\partial_{xx}\psi_k( \check{\eta}^{k}_{s_{k}})\end{pmatrix}
 =\begin{pmatrix} 
  b_2,& 2\beta^{\frac{1}{3}}(\hat{{\eta}}_{{\hat{t}}}({{\hat{t}}})-\hat{{\gamma}}_{{\hat{t}}}({{\hat{t}}})),& Y\end{pmatrix},
\end{eqnarray}
 where $b_{1}+b_{2}=0$ and $X,Y\in \mathcal{S}(\mathbb{R}^{d})$ satisfy the following inequality:
\begin{eqnarray}\label{II}
                              {-6\beta^{\frac{1}{3}}}\left(\begin{array}{cc}
                                    I&0\\
                                    0&I
                                    \end{array}\right)\leq \left(\begin{array}{cc}
                                    X&0\\
                                    0&Y
                                    \end{array}\right)\leq  6\beta^{\frac{1}{3}} \left(\begin{array}{cc}
                                    I&-I\\
                                    -I&I
                                    \end{array}\right).
\end{eqnarray}
            We note that   sequence  $(\check{\gamma}^{k}_{l_k},\check{\eta}^{k}_{s_k},l_{k},s_{k},\varphi_k,\psi_k)$ and $b_{1},b_{2},X,Y$ depend on  $\beta$ and
              $\varepsilon$.
\par
 For every $(t,\gamma_t),(s,\eta_s)\in [T-\bar{a},T]\times{{\Lambda}^{T-\bar{a}}}$, let
\begin{eqnarray*}
\chi^{k}(\gamma_t)
        :=
                \varepsilon\frac{\nu T-t}{\nu
                 T}\Upsilon^{m}(\gamma_t)
                +\varepsilon \overline{\Upsilon}_0(\gamma_t,\hat{\gamma}_{\hat{t}})
                +\sum_{i=0}^{\infty}
        \frac{1}{2^i}\overline{\Upsilon}_0(\gamma^i_{t_i},\gamma_t)+2^5\beta\Upsilon_0(\gamma_{t},\hat{{\xi}}_{{\hat{t}}})
                +\varphi_k(\gamma_t),
  \end{eqnarray*}
  \begin{eqnarray*}
\hbar^{k}(\eta_s)
        :=-\varepsilon\frac{\nu T-s}{\nu
                 T}\Upsilon^{m}(\eta_s)
                 -\varepsilon \overline{\Upsilon}_0(\eta_s,\hat{\eta}_{\hat{t}})
                 -\sum_{i=0}^{\infty}
        \frac{1}{2^i}{\Upsilon}_0(\eta^i_{t_i},\eta_s)-2^5\beta\Upsilon_0(\eta_{s},\hat{{\xi}}_{{\hat{t}}})
        -\psi_k(\eta_s).
  \end{eqnarray*}
  Then $\chi^{k}(\cdot)\in C^{1,2}_p(\Lambda^{l_k}),\hbar^{k}(\cdot)\in C^{1,2}_p(\Lambda^{s_k})$. 
  Moreover, by  (\ref{0609a}), (\ref{0609b}) and definitions of $w_1$ and $w_2$,
  $$
                         (W_1-\chi^{k})(\check{\gamma}^{k}_{l_k})=\sup_{(t,\gamma_t)\in [l_k,T]\times\Lambda^{l_k}}
                         (W_1-\chi^{k})(\gamma_t),
$$
$$
                         (W_2-\hbar^{k})(\check{\eta}^{k}_{s_k})=\inf_{(s,\eta_s)\in [s_k,T]\times\Lambda^{s_k}}
                         (W_2-\hbar^{k})(\eta_s).
$$
From $l_{k}\rightarrow {\hat{t}}$, $s_{k}\rightarrow {\hat{t}}$ as $k\rightarrow\infty$ and ${\hat{t}}<T$
 for $\beta\geq N$, it follows that, for every fixed $\beta\geq N$,   constant $ K_\beta>0$ exists such that
$$
            l_{k}\vee s_{k}<T, 
            \ \ \mbox{for all}    \ \ k\geq K_\beta.
$$
Now, for every $\beta\geq N$ and $k>K_\beta$, from the definition of viscosity solutions it follows that
  \begin{eqnarray}\label{vis1}
                      \partial_t\chi^{k}(\check{\gamma}^{k}_{l_{k}})
  +{\mathbf{F}}{(}\check{\gamma}^{k}_{l_{k}},  W_1(\check{\gamma}^{k}_{l_{k}}), \partial_x\chi^{k}(\check{\gamma}^{k}_{l_{k}}),
                                       \partial_{xx}\chi^{k}(\check{\gamma}^{k}_{l_{k}})
                                    {)}\geq c,
 \end{eqnarray}
 and
  \begin{eqnarray}\label{vis2}
                     \partial_t\hbar^{k}(\check{\eta}^{k}_{s_{k}})+{\mathbf{F}}{(}\check{\eta}^{k}_{s_{k}}, W_2(\check{\eta}^{k}_{s_{k}}),
                     \partial_x\hbar^{k}(\check{\eta}^{k}_{s_{k}}),\partial_{xx}\hbar^{k}(\check{\eta}^{k}_{s_{k}}){)}\leq0,
  \end{eqnarray}
 where, for every $(t,\gamma_t)\in [l_k,T]\times{{\Lambda}^{l_k}}$ and  $ (s,\eta_s)\in [s_k,T]\times{{\Lambda}^{s_k}}$,  from Lemma \ref{theoremS},
  \begin{eqnarray*}
\partial_t\chi^{k}(\gamma_t)=
                       \partial_t\varphi_k(\gamma_{{t}})-\frac{\varepsilon}{\nu T}\Upsilon^{m}(\gamma_{{t}})
                     +2\varepsilon({t}-{\hat{t}})+2\sum_{i=0}^{\infty}\frac{1}{2^i}(t-t_{i}),
  \end{eqnarray*}
  \begin{eqnarray*}
\partial_x\chi^{k}(\gamma_t)&=&\partial_{x}\varphi_k(\gamma_{{t}})
                     +\varepsilon\frac{\nu T-{t}}{\nu T}\partial_x\Upsilon^{m}(\gamma_{{t}}) +\varepsilon\partial_x\Upsilon(\gamma_{{t}}-\hat{\gamma}_{{\hat{t},t}})
+2^5\beta\partial_x\Upsilon(\gamma_{{t}}-\hat{\xi}_{{\hat{t},t}})
\\
 &&+\partial_x\left[\sum_{i=0}^{\infty}\frac{1}{2^i}
                      \Upsilon(\gamma_{t}-\gamma^{i}_{t_{i},t})
                     \right],
  \end{eqnarray*}
   \begin{eqnarray*}
\partial_{xx}\chi^{k}(\gamma_t)&=&\partial_{xx}(\varphi_k)(\gamma_{{t}})+\varepsilon\frac{\nu T-{t}}{\nu T}\partial_{xx}\Upsilon^{m}(\gamma_{{t}}) +\varepsilon\partial_{xx}\Upsilon(\gamma_{{t}}-\hat{\gamma}_{{\hat{t},t}})
+2^5\beta\partial_{xx}\Upsilon(\gamma_{{t}}-\hat{\xi}_{{\hat{t},t}})
\\
                      &&
                      +\partial_{xx}\left[\sum_{i=0}^{\infty}\frac{1}{2^i}
                      \Upsilon(\gamma_{t}-\gamma^{i}_{t_{i},t})
                      \right],
  \end{eqnarray*}
   \begin{eqnarray*}
\partial_t\hbar^{k}(\eta_s)=
                      -\partial_t\psi_k(\eta_{{s}})+\frac{\varepsilon}{\nu T}\Upsilon^{m}(\eta_{{s}})-2\varepsilon({s}-{\hat{t}}),
  \end{eqnarray*}
  \begin{eqnarray*}
\partial_x\hbar^{k}(\eta_s)
                          &=& -\partial_{x}\psi_k(\eta_{{s}})
                      -\varepsilon\frac{\nu T-{s}}{\nu T}
                     \partial_x\Upsilon^{m}(\eta_{{s}})-\varepsilon\partial_x\Upsilon(\eta_{{s}}-\hat{\eta}_{{\hat{t}},s})
                     -2^5\beta\partial_x\Upsilon(\eta_{{s}}-\hat{\xi}_{{\hat{t}},s})\\
                     &&
                      -\partial_x\left[\sum_{i=0}^{\infty}\frac{1}{2^i}
                      \Upsilon(\eta_{s}-{\eta}^{i}_{{t}_{i},s})
                      \right],
  \end{eqnarray*}
   \begin{eqnarray*}
\partial_{xx}\hbar^{k}(\eta_s)
                                &=&-\partial_{xx}\psi_k( \eta_{{s}})-\varepsilon\frac{\nu T-{s}}{\nu T}
                     \partial_{xx}\Upsilon^{m}(\eta_{{s}})-\varepsilon\partial_{xx}\Upsilon(\eta_{{s}}-\hat{\eta}_{{\hat{t}},s})-2^5\beta\partial_{xx}\Upsilon(\eta_{{s}}-\hat{\xi}_{{\hat{t}},s})\\
                     &&
                      -\partial_{xx}\left[\sum_{i=0}^{\infty}\frac{1}{2^i}
                      \Upsilon(\eta_{s}-{\eta}^{i}_{{t}_{i},s})
                      \right].
  \end{eqnarray*}
 \par

    Letting  $k\rightarrow\infty$ in (\ref{vis1}) and (\ref{vis2}), and using (\ref{s0}), Hypothesis  \ref{hypstate} (i), (\ref{4.23}), (\ref{0608v1}) and  (\ref{0608vw1}),    we obtain
    \begin{eqnarray}\label{03103}
                   b_1-\frac{\varepsilon}{\nu T}\Upsilon^{m}(\hat{{\gamma}}_{{\hat{t}}})+2\sum_{i=0}^{\infty}\frac{1}{2^i}(\hat{t}-t_i)
                     +{\mathbf{F}}(\hat{{\gamma}}_{{\hat{t}}},W_1(\hat{{\gamma}}_{{\hat{t}}}),
                    \partial_x\chi(\hat{\gamma}_{\hat{t}}),\partial_{xx}\chi(\hat{\gamma}_{\hat{t}}))
                     \geq c;
\end{eqnarray}
and
 \begin{eqnarray}\label{03104}
                     -b_2+ \frac{\varepsilon}{\nu T}\Upsilon^{m}(\hat{{\eta}}_{{\hat{t}}})+{\mathbf{F}}(\hat{{\eta}}_{{\hat{t}}},W_2(\hat{{\eta}}_{{\hat{t}}}),\partial_x\hbar(\hat{\eta}_{{\hat{t}}}),\partial_{xx}\hbar(\hat{\eta}_{\hat{t}}))
                     \leq0,
\end{eqnarray}
  where
  \begin{eqnarray*}
\partial_x\chi(\hat{\gamma}_{\hat{t}})
                                     :=
2\beta^{\frac{1}{3}}(\hat{{\gamma}}_{{\hat{t}}}({{\hat{t}}})-\hat{{\eta}}_{{\hat{t}}}({{\hat{t}}}))+2^5\beta\partial_x\Upsilon(\hat{\gamma}_{{\hat{t}}}-\hat{\xi}_{{\hat{t}}})
                                      +\varepsilon\frac{\nu T-{\hat{t}}}{\nu T}\partial_x\Upsilon^{m}(\hat{{\gamma}}_{{\hat{t}}})
                                    +\partial_x\left[\sum_{i=0}^{\infty}\frac{1}{2^i}\Upsilon(\hat{\gamma}_{\hat{t}}-\gamma^i_{t_i,\hat{t}})\right],
  \end{eqnarray*}
   \begin{eqnarray*}
\partial_{xx}\chi(\hat{\gamma}_{\hat{t}})
:= X+2^5\beta\partial_{xx}\Upsilon(\hat{\gamma}_{{\hat{t}}}-\hat{\xi}_{{\hat{t}}})
                                     +\varepsilon\frac{\nu T-{\hat{t}}}{\nu T}\partial_{xx}\Upsilon^{m}(\hat{{\gamma}}_{{\hat{t}}})
                                     +\partial_{xx}\left[\sum_{i=0}^{\infty}\frac{1}{2^i}\Upsilon(\hat{\gamma}_{\hat{t}}-\gamma^i_{t_i,\hat{t}})\right],
  \end{eqnarray*}
  \begin{eqnarray*}
\partial_x\hbar(\hat{\eta}_{{\hat{t}}})
                                     :=
2\beta^{\frac{1}{3}}(\hat{{\gamma}}_{{\hat{t}}}({{\hat{t}}})-\hat{{\eta}}_{{\hat{t}}}({{\hat{t}}}))-2^5\beta\partial_x\Upsilon(\hat{\eta}_{{\hat{t}}}-\hat{\xi}_{{\hat{t}}})-\varepsilon\frac{\nu T-{\hat{t}}}{\nu T}\partial_x\Upsilon^{m}(\hat{{\eta}}_{{\hat{t}}})-\partial_x\left[\sum_{i=0}^{\infty}\frac{1}{2^i}
                                     \Upsilon(\hat{\eta}_{\hat{t}}-\eta^i_{t_i,\hat{t}})\right],
  \end{eqnarray*}
  and
   \begin{eqnarray*}
\partial_{xx}\hbar(\hat{\eta}_{\hat{t}})
:= -Y-2^5\beta\partial_{xx}\Upsilon(\hat{\eta}_{{\hat{t}}}-\hat{\xi}_{{\hat{t}}})-\varepsilon\frac{\nu T-{\hat{t}}}{\nu T}\partial_{xx}\Upsilon^{m}(\hat{{\eta}}_{{\hat{t}}})
                      -\partial_{xx}\left[\sum_{i=0}^{\infty}\frac{1}{2^i}
                      \Upsilon(\hat{\eta}_{\hat{t}}-\eta^i_{t_i,\hat{t}})\right].
  \end{eqnarray*}
  Notice that $b_1+b_2=0$ and $\hat{\xi}_{{\hat{t}}}=\frac{\hat{{\gamma}}_{\hat{t}}+\hat{{\eta}}_{\hat{t}}}{2}$,
combining  (\ref{03103}) and (\ref{03104}),  we have
 \begin{eqnarray}\label{vis112}
                    && c+ \frac{\varepsilon}{\nu T}(\Upsilon^{m}(\hat{{\gamma}}_{{\hat{t}}})+\Upsilon^{m}(\hat{{\eta}}_{{\hat{t}}})
                     )-2\sum_{i=0}^{\infty}\frac{1}{2^i}(\hat{t}-t_i)\nonumber\\
                     &\leq&{\mathbf{F}}(\hat{{\gamma}}_{{\hat{t}}},W_1(\hat{{\gamma}}_{{\hat{t}}}),\partial_x\chi(\hat{\gamma}_{\hat{t}}),\partial_{xx}\chi(\hat{\gamma}_{\hat{t}}))
                     -{\mathbf{F}}(\hat{{\eta}}_{{\hat{t}}},W_2(\hat{{\eta}}_{{\hat{t}}}),\partial_x\hbar(\hat{\eta}_{{\hat{t}}}),\partial_{xx}\hbar(\hat{\eta}_{\hat{t}})).
\end{eqnarray}
 On the other hand, by  Hypothesis \ref{hypstate} (ii) and via a simple calculation we obtain
 \begin{eqnarray}\label{v4}
                &&{\mathbf{F}}(\hat{{\gamma}}_{{\hat{t}}},W_1(\hat{{\gamma}}_{{\hat{t}}}),\partial_x\chi(\hat{\gamma}_{\hat{t}}),\partial_{xx}\chi(\hat{\gamma}_{\hat{t}}))
                     -{\mathbf{F}}(\hat{{\eta}}_{{\hat{t}}},W_2(\hat{{\eta}}_{{\hat{t}}}),\partial_x\hbar(\hat{\eta}_{{\hat{t}}}),\partial_{xx}\hbar(\hat{\eta}_{\hat{t}}))\nonumber\\
                     &\leq&{\mathbf{F}}(\hat{{\gamma}}_{{\hat{t}}},W_2(\hat{{\eta}}_{{\hat{t}}}),\partial_x\chi(\hat{\gamma}_{\hat{t}}),\partial_{xx}\chi(\hat{\gamma}_{\hat{t}}))
                     -{\mathbf{F}}(\hat{{\eta}}_{{\hat{t}}},W_2(\hat{{\eta}}_{{\hat{t}}}),\partial_x\hbar(\hat{\eta}_{{\hat{t}}}),\partial_{xx}\hbar(\hat{\eta}_{\hat{t}}))\nonumber\\
                &=& {\mathbf{F}}(\hat{{\gamma}}_{{\hat{t}}},W_2(\hat{{\eta}}_{{\hat{t}}}),\partial_x\chi(\hat{\gamma}_{\hat{t}}),\partial_{xx}\chi(\hat{\gamma}_{\hat{t}}))
                -{\mathbf{F}}(\hat{{\gamma}}_{{\hat{t}}},W_2(\hat{{\eta}}_{{\hat{t}}}),2\beta^{\frac{1}{3}}(\hat{{\gamma}}_{{\hat{t}}}({{\hat{t}}})-\hat{{\eta}}_{{\hat{t}}}({{\hat{t}}})),X)\nonumber\\
                     &&+{\mathbf{F}}(\hat{{\gamma}}_{{\hat{t}}},W_2(\hat{{\eta}}_{{\hat{t}}}),2\beta^{\frac{1}{3}}(\hat{{\gamma}}_{{\hat{t}}}({{\hat{t}}})-\hat{{\eta}}_{{\hat{t}}}({{\hat{t}}})),X)
                     -{\mathbf{F}}(\hat{{\eta}}_{{\hat{t}}},W_2(\hat{{\eta}}_{{\hat{t}}}),2\beta^{\frac{1}{3}}(\hat{{\gamma}}_{{\hat{t}}}({{\hat{t}}})-\hat{{\eta}}_{{\hat{t}}}({{\hat{t}}})),-Y)\nonumber\\
                     &&+{\mathbf{F}}(\hat{{\eta}}_{{\hat{t}}},W_2(\hat{{\eta}}_{{\hat{t}}}),2\beta^{\frac{1}{3}}(\hat{{\gamma}}_{{\hat{t}}}({{\hat{t}}})-\hat{{\eta}}_{{\hat{t}}}({{\hat{t}}})),-Y)
                     -{\mathbf{F}}(\hat{{\eta}}_{{\hat{t}}},W_2(\hat{{\eta}}_{{\hat{t}}}),\partial_x\hbar(\hat{\eta}_{{\hat{t}}}),\partial_{xx}\hbar(\hat{\eta}_{\hat{t}}))\nonumber\\
                     &=& J_{1}+J_{2}+J_{3},
\end{eqnarray}
 where, from Hypothesis \ref{hypstate} (v) and  (\ref{s0}), 
\begin{eqnarray}\label{j1}
                J_1
                     &\leq& M_F(1+||\hat{{\gamma}}_{{\hat{t}}}||_0)\left|2^5\beta\partial_x\Upsilon(\hat{\gamma}_{{\hat{t}}}-\hat{\xi}_{{\hat{t}}})
                                      +\varepsilon\frac{\nu T-{\hat{t}}}{\nu T}\partial_x\Upsilon^{m}(\hat{{\gamma}}_{{\hat{t}}})
                                    +\partial_x\left[\sum_{i=0}^{\infty}\frac{1}{2^i}\Upsilon(\hat{\gamma}_{\hat{t}}-\gamma^i_{t_i,\hat{t}})\right]\right|\nonumber\\
                                    &&+M_F(1+||\hat{{\gamma}}_{{\hat{t}}}||_0^2)\left|2^5\beta\partial_{xx}\Upsilon(\hat{\gamma}_{{\hat{t}}}-\hat{\xi}_{{\hat{t}}})
                                     +\varepsilon\frac{\nu T-{\hat{t}}}{\nu T}\partial_{xx}\Upsilon^{m}(\hat{{\gamma}}_{{\hat{t}}})
                                     +\partial_{xx}\left[\sum_{i=0}^{\infty}\frac{1}{2^i}\Upsilon(\hat{\gamma}_{\hat{t}}-\gamma^i_{t_i,\hat{t}})\right]\right|\nonumber\\
                      &\leq&18M_F(1+||\hat{{\gamma}}_{{\hat{t}}}||_0)\left(\beta|\hat{{\gamma}}_{{\hat{t}}}({\hat{t}})-\hat{{\eta}}_{{\hat{t}}}({\hat{t}})|^5+
                      \sum_{i=0}^{\infty}\frac{1}{2^i}|\gamma^i_{t_i}(t_i)-\hat{\gamma}_{\hat{t}}(\hat{t})|^5
                                            \right)\nonumber\\
                                    &&+306M_F(1+||\hat{{\gamma}}_{{\hat{t}}}||_0^2)\left(2\beta|\hat{{\gamma}}_{{\hat{t}}}({\hat{t}})-\hat{{\eta}}_{{\hat{t}}}({\hat{t}})|^4+
                                     \sum_{i=0}^{\infty}\frac{1}{2^i}\left|\gamma^i_{t_i}(t_i)-\hat{\gamma}_{\hat{t}}(\hat{t})\right|^4
                                     \right)\nonumber\\
                                     &&+6m\varepsilon \frac{\nu T-{\hat{t}}}{\nu T} M_F(1+2||\hat{{\gamma}}_{{\hat{t}}}||^{2m}_0)+6m(6m-1)\varepsilon \frac{\nu T-{\hat{t}}}{\nu T}M_F(1+2||\hat{{\gamma}}_{{\hat{t}}}||^{2m}_0);
\end{eqnarray}
\begin{eqnarray}\label{j3}
                J_3
                     &\leq&
                               M_F(1+||\hat{{\eta}}_{{\hat{t}}}||_0)\left|2^5\beta\partial_x\Upsilon(\hat{\eta}_{{\hat{t}}}-\hat{\xi}_{{\hat{t}}})+\varepsilon\frac{\nu T-{\hat{t}}}{\nu T}\partial_x\Upsilon^{m}(\hat{{\eta}}_{{\hat{t}}})+\partial_x\left[\sum_{i=0}^{\infty}\frac{1}{2^i}
                                     \Upsilon(\hat{\eta}_{\hat{t}}-\eta^i_{t_i,\hat{t}})\right]\right|\nonumber\\
                                    &&+M_F(1+||\hat{{\eta}}_{{\hat{t}}}||_0^2)\left|2^5\beta\partial_{xx}\Upsilon(\hat{\eta}_{{\hat{t}}}-\hat{\xi}_{{\hat{t}}})+\varepsilon\frac{\nu T-{\hat{t}}}{\nu T}\partial_{xx}\Upsilon^{m}(\hat{{\eta}}_{{\hat{t}}})
                      +\partial_{xx}\left[\sum_{i=0}^{\infty}\frac{1}{2^i}
                      \Upsilon(\hat{\eta}_{\hat{t}}-\eta^i_{t_i,\hat{t}})\right]\right|\nonumber\\
                      &\leq&18M_F(1+||\hat{{\eta}}_{{\hat{t}}}||_0)\left(\beta|\hat{{\gamma}}_{{\hat{t}}}({\hat{t}})-\hat{{\eta}}_{{\hat{t}}}({\hat{t}})|^5+
                      \sum_{i=0}^{\infty}\frac{1}{2^i}|\eta^i_{t_i}(t_i)-\hat{\eta}_{\hat{t}}(\hat{t})|^5\right)\nonumber\\
                                    &&+306M_F(1+||\hat{{\eta}}_{{\hat{t}}}||_0^2)\left(2\beta|\hat{{\gamma}}_{{\hat{t}}}({\hat{t}})-\hat{{\eta}}_{{\hat{t}}}({\hat{t}})|^4+
                                     \sum_{i=0}^{\infty}\frac{1}{2^i}\left|\eta^i_{t_i}(t_i)-\hat{\eta}_{\hat{t}}(\hat{t})\right|^4\right)\nonumber\\
                                     &&+6m\varepsilon \frac{\nu T-{\hat{t}}}{\nu T} M_F(1+2||\hat{{\eta}}_{{\hat{t}}}||^{2m}_0)+6m(6m-1)\varepsilon \frac{\nu T-{\hat{t}}}{\nu T}M_F(1+2||\hat{{\eta}}_{{\hat{t}}}||^{2m}_0);
\end{eqnarray}
and, from Hypothesis \ref{hypstate} (iv), (\ref{w}), (\ref{5.10jiajiaaaa}) and (\ref{II})
\begin{eqnarray}\label{j2}
                J_2
                     &\leq& \rho(2
                               \beta^{\frac{1}{3}}||\hat{{\gamma}}_{{\hat{t}}}-\hat{{\eta}}_{{\hat{t}}}||_0^2+||\hat{{\gamma}}_{{\hat{t}}}-\hat{{\eta}}_{{\hat{t}}}||_0,
                               |W_2(\hat{{\eta}}_{{\hat{t}}})|\vee||\hat{{\gamma}}_{{\hat{t}}}||_0\vee||\hat{{\eta}}_{{\hat{t}}}||_0)\nonumber\\
                               &\leq& \rho(2
                               \beta^{\frac{1}{3}}||\hat{{\gamma}}_{{\hat{t}}}-\hat{{\eta}}_{{\hat{t}}}||_0^2+||\hat{{\gamma}}_{{\hat{t}}}-\hat{{\eta}}_{{\hat{t}}}||_0,L(1+M_0^m)\vee M_0).
\end{eqnarray}
   We have from  the property (i) of $(\hat{t},(\hat{\gamma}_{\hat{t}},\hat{\eta}_{\hat{t}}))$ that
  \begin{eqnarray*}
2\sum_{i=0}^{\infty}\frac{1}{2^i}(\hat{t}-t_i)
  \leq2\sum_{i=0}^{\infty}\frac{1}{2^i}\bigg{(}\frac{1}{2^i\beta}\bigg{)}^{\frac{1}{2}}\leq 4{\bigg{(}\frac{1}{{\beta}}\bigg{)}}^{\frac{1}{2}},
    \end{eqnarray*}
     \begin{eqnarray*}
  \sum_{i=0}^{\infty}\frac{1}{2^i}\left[\left|\gamma^i_{t_i}(t_i)-\hat{\gamma}_{\hat{t}}(\hat{t})\right|^5
                                            +\left|\eta^i_{t_i}(t_i)-\hat{\eta}_{\hat{t}}(\hat{t})\right|^5\right]
                      \leq
                      2\sum_{i=0}^{\infty}\frac{1}{2^i}\bigg{(}\frac{1}{2^i\beta}\bigg{)}^{\frac{5}{6}}\leq 4{\bigg{(}\frac{1}{{\beta}}\bigg{)}}^{\frac{5}{6}},
                        \end{eqnarray*}
                        and
                         \begin{eqnarray*}
  \sum_{i=0}^{\infty}\frac{1}{2^i}\left[\left|\gamma^i_{t_i}(t_i)-\hat{\gamma}_{\hat{t}}(\hat{t})\right|^4
                                     +\left|\eta^i_{t_i}(t_i)-\hat{\eta}_{\hat{t}}(\hat{t})\right|^4\right]
                      \leq
                      2\sum_{i=0}^{\infty}\frac{1}{2^i}\bigg{(}\frac{1}{2^i\beta}\bigg{)}^{\frac{2}{3}}\leq 4{\bigg{(}\frac{1}{{\beta}}\bigg{)}}^{\frac{2}{3}}.
  \end{eqnarray*}
 Combining (\ref{vis112})-(\ref{j2}),  we see from  (\ref{5.10jiajiaaaa}) and (\ref{5.10})  that for sufficiently large  $\beta>0$,
\begin{eqnarray}\label{vis122}
                     c
                                             \leq
                      -\frac{\varepsilon}{\nu T}(\Upsilon^{m}(\hat{{\gamma}}_{{\hat{t}}})
                     +\Upsilon^{m}(\hat{{\eta}}_{{\hat{t}}}))+ 72m^2\varepsilon \frac{\nu T-{\hat{t}}}{\nu T} M_F(1+||\hat{{\gamma}}_{{\hat{t}}}||^{2m}_0+||\hat{{\eta}}_{{\hat{t}}}||^{2m}_0)+\frac{c}{4}.
\end{eqnarray}
Since $
         \nu=1+\frac{1}{144m^2M_FT}
$ and $\bar{a}=\frac{1}{144m^2M_F}\wedge{T}$,  (\ref{up}) and (\ref{5.3}) lead to the following contradiction:
\begin{eqnarray*}\label{vis122}
                     c\leq
                             \frac{\varepsilon}{\nu
                              T}+\frac{c}{4}\leq \frac{c}{2}.
\end{eqnarray*}
 The proof is complete.
 \ \ $\Box$

\section{Application to path-dependent stochastic differential games.}
\par
Let $\Omega:=\{\omega\in C([0,T],\mathbb{R}^n):\omega(0)={\mathbf{0}}\}$, the set of continuous functions with initial value ${\mathbf{0}}$,
           $W$ the canonical process, $\mathbb{P}$ the Wiener measure, ${\cal {F}}$ the  Borel $\sigma$-field over $\Omega$, completed with
            respect to the Wiener measure $\mathbb{P}$ on this space.  Then  $(\Omega,{\cal {F}},\mathbb{P})$ is a complete   space.   By $\mathbb{F}=\{{\cal {F}}_t\}_{0\leq t\leq T}$ we denote  the filtration generated by $\{W(t),0\leq t\leq T\}$, augmented
       with the family $\mathcal {N}$ of $\mathbb{P}$-null of ${\cal {F}}$.
 The filtration $\mathbb{F}$ satisfies the
       usual conditions.
       \par
We introduce the admissible control and admissible strategy. Let $t,s$ be two deterministic times, $0\leq t<s\leq T$.
\begin{definition}
An admissible control process $u=\{u(r):  r\in [t,s]\}$ (resp., $v=\{v(r):  r\in [t,s]\}$) for player I (resp., II) on $[t,s]$  ($t<s \leq T$) is an ${\mathbb{F}}$-progressively
measurable process taking values in some $\sigma$-compact metric space $(U,d_1)$ (resp., $(V,d_2)$). The set of all admissible controls for
player I (resp., II) on $[t,s]$ is denoted by ${\cal{U}}[t,s]$ (resp., ${\cal{V}}[t,s]$). We identify two processes $u$ and $\tilde{u}$ in ${\cal{U}}[t,s]$
                and write $u\equiv\tilde{u}$ on $[t,s]$, if $\mathbb{P}(u=\tilde{u} \ a.e. \ \mbox{in}\ [t,s])=1$. Similarly we
interpret $v\equiv\tilde{v}$ on $[t,s]$ in ${\cal{V}}[t,s]$.
\end{definition}
\begin{definition}
                A nonanticipative strategy for player I on $[t, s]$ ($t < s \leq T$)
is a mapping $\alpha : {\cal{V}}[t,s] \rightarrow {\cal{U}}[t,s]$ such that, for any ${\mathbb{F}}$-stopping time $S : \Omega\rightarrow[t,s]$
and any $v_1, v_2 \in {\cal{V}}[t,s]$, with $v_1\equiv v_2$ on $[[t, S]]$, it holds that $\alpha(v_1) \equiv\alpha(v_2)$ on $[[t, S]]$.
Nonanticipative strategies for player II on $[t, s]$, $\beta : {\cal{U}}[t,s] \rightarrow {\cal{V}}[t,s]$, are defined similarly.
The set of all nonanticipative strategies $\alpha$ (resp., $\beta$) for player I  (resp., II) on $[t, s]$ is denoted
by ${\cal{A}}_{[t,s]}$ (resp., ${\cal{B}}_{[t,s]}$).
\end{definition}
\par
We consider the following   controlled path-dependent stochastic  differential
                 game (PSDG):
\begin{eqnarray}\label{state1}
\begin{cases}
            dX^{\gamma_t,u,v}(s)=
           b(X_s^{\gamma_t,u,v},u(s),v(s))ds+\sigma(X_s^{\gamma_t,u,v},u(s),v(s))dW(s),  \ \ s\in [t,T],\\
~~~~~X_t^{\gamma_t,u,v}=\gamma_t\in {\Lambda}_t.
\end{cases}
\end{eqnarray}
The cost functional (interpreted as a payoff for player I and as a cost for player
II) is introduced by a backward stochastic differential equation (BSDE):
\begin{eqnarray}\label{fbsde1}
Y^{\gamma_t,u,v}(s)&=&\phi(X_T^{\gamma_t,u,v})+\int^{T}_{s}q(X_\sigma^{\gamma_t,u,v},Y^{\gamma_t,u,v}(\sigma),Z^{\gamma_t,u,v}(\sigma),u(\sigma),v(\sigma))d\sigma\nonumber\\
                 &&-\int^{T}_{s}Z^{\gamma_t,u,v}(\sigma)dW(\sigma),\ \ \ a.s., \ \ \mbox{all}\ \ s\in [t,T].
\end{eqnarray}
The payoff is given by
\begin{eqnarray}\label{cost1}
                     J(\gamma_t,u,v):=Y^{\gamma_t,u,v}(t),\ \ \ (t,\gamma_t)\in [0,T]\times {\Lambda}.
\end{eqnarray}
 We define the lower value functional of our PSDG:
\begin{eqnarray}\label{value1}
\underline{V}(\gamma_t):=\mathop{\essinf}\limits_{\beta(\cdot)\in {\cal{B}}[t,T]}\mathop{\esssup}\limits_{u(\cdot)\in{\cal{U}}[t,T]} Y^{\gamma_t,u,\beta(u)}(t),\ \  (t,\gamma_t)\in [0,T]\times {\Lambda},
\end{eqnarray}
and its upper value functional:
\begin{eqnarray}\label{value2}
\overline{V}(\gamma_t):=\mathop{\esssup}\limits_{\alpha(\cdot)\in {\cal{A}}[t,T]}\mathop{\essinf}\limits_{v(\cdot)\in{\cal{V}}[t,T]} Y^{\gamma_t,\alpha(v),v}(t),\ \  (t,\gamma_t)\in [0,T]\times {\Lambda}.
\end{eqnarray}

For $(t,\gamma_t,r,p,l)\in [0,T]\times{\Lambda}\times \mathbb{R}\times \mathbb{R}^d\times {\cal{S}}(\mathbb{R}^{d})$, define
\begin{eqnarray*}
	{\mathbf{H}}^-(\gamma_t,r,p,l):=\sup_{u\in U}\inf_{v\in V} \mathbf{H}(\gamma_t,r,p,l,u,v), \\
	{\mathbf{H}}^+(\gamma_t,r,p,l):=\inf_{v\in V}\sup_{u\in
		U} \mathbf{H}(\gamma_t,r,p,l,u,v),
\end{eqnarray*}
and
\begin{eqnarray*}
	\mathbf{H}(\gamma_t,r,p,l,u,v)&:=&               [
	(p,b(\gamma_t,u,v))_{\mathbb{R}^d}+\frac{1}{2}\mbox{tr}[ l \sigma(\gamma_t,u,v)\sigma^\top(\gamma_t,u,v)]\\
	&&\ \ \    +q(\gamma_t,r,\sigma^\top(\gamma_t,u,v)p,u,v)].
\end{eqnarray*}

  The goal of this section  is to characterize the lower  and  upper value functionals $\underline{V}$ and $\overline{V}$ as the unique viscosity solution of the following
                    path-dependent  Hamilton-Jacobi-Bellman-Isaacs equations (PHJBIEs):
  \begin{eqnarray}\label{hjb106}
\begin{cases}
\partial_t\underline{V}(\gamma_t)+{\mathbf{H}}^-(\gamma_t,\underline{V}(\gamma_t),\partial_x\underline{V}(\gamma_t),\partial_{xx}\underline{V}(\gamma_t))= 0,\ \ \  (t,\gamma_t)\in
                               [0,T)\times {\Lambda};\\
 \underline{V}(\gamma_T)=\phi(\gamma_T), \ \ \ \gamma_T\in {\Lambda}_T
 \end{cases}
\end{eqnarray}
and
\begin{eqnarray}\label{hjb2}
\begin{cases}
\partial_t\overline{V}(\gamma_t)+{\mathbf{H}}^+(\gamma_t,\overline{V}(\gamma_t),\partial_x\overline{V}(\gamma_t),\partial_{xx}\overline{V}(\gamma_t))= 0,\ \ \  (t,\gamma_t)\in
                               [0,T)\times {\Lambda};\\
 \overline{V}(\gamma_T)=\phi(\gamma_T), \ \ \ \gamma_T\in {\Lambda}_T,
 \end{cases}
\end{eqnarray}
respectively.

We make the following assumption.
\begin{hyp}\label{hypstate6}
$b:{\Lambda}\times U\times V\rightarrow \mathbb{R}^{d}$, $\sigma:{\Lambda}\times U\times V\rightarrow \mathbb{R}^{d\times n}$, $
        q: {\Lambda}\times \mathbb{R}\times \mathbb{R}^n\times U\times V\rightarrow \mathbb{R}$ and $\phi: {\Lambda}_T\rightarrow \mathbb{{R}}$ are continuous, and $b, \sigma, q$ are continuous in $\gamma_t\in \Lambda$, uniformly in $(u,v)\in U\times V$. Moreover,
                 there exists   $L>0$
                 such that, for all $(t,\gamma_t,\eta_T,y,z,u,v)$,  $ (t, \gamma'_t,\eta'_T,y',z',u,v)
                 \in [0,T]\times\Lambda\times {\Lambda}_T\times \mathbb{R}\times \mathbb{R}^n\times U\times V$,
      \begin{eqnarray}\label{211224}
                |b(\gamma_t,u,v)|^2\vee|\sigma(\gamma_t,u,v)|_2^2\leq
                 L^2(1+||\gamma_t||_0^2),
        \end{eqnarray}
         \begin{eqnarray}\label{2112241}
                 |b(\gamma_t,u,v)-b(\gamma'_{t},u,v)|\vee|\sigma(\gamma_t,u,v)-\sigma(\gamma'_{t},u,v)|_2\leq
                 L||\gamma_t-\gamma'_{t}||_0,
            \end{eqnarray}
             \begin{eqnarray}\label{2112242}
                  |q(\gamma_t,y,z,u,v)|\leq L(1+||\gamma_t||_0+|y|+|z|),
                 \end{eqnarray}
                  \begin{eqnarray}\label{2112243}
               |q(\gamma_t,y,z,u,v)-q(\gamma'_{t},y',z',u,v)|\leq L(||\gamma_t-\gamma'_{t}||_0+|y-y'|+|z-z'|),
                \end{eqnarray}
                  \begin{eqnarray}\label{2112244}
                 |\phi(\eta_T)-\phi(\eta'_T)|\leq L||\eta_T-\eta'_{T}||_0.
\end{eqnarray}
\end{hyp}
\begin{lemma}\label{lemmaexist} (\cite[Lemma 2.3]{zhang})
\ \ Assume that Hypothesis \ref{hypstate6}  holds. Then for every $u(\cdot)\in {\cal{U}}[0,T]$, $v(\cdot)\in {\cal{V}}[0,T]$ ,
$(t,\gamma_t)\in [0,T]\times {\Lambda}$ and $p\geq2$, PSDE (\ref{state1}) admits a
unique strong  solution $X^{\gamma_t,u,v}$, and BSDE (\ref{fbsde1}) admits a unique
 pair of solutions $(Y^{\gamma_t,u,v}, Z^{\gamma_t,u,v})$.  Furthermore, let  $X^{\gamma'_t,u,v}$ and $(Y^{\gamma'_t,u,v}, Z^{\gamma'_t,u,v})$ be the solutions of PSDE (\ref{state1}) and BSDE (\ref{fbsde1})
 corresponding $(t,\gamma'_t)\in [0,T]\times \Lambda$, $u(\cdot)\in {\cal{U}}[0,T]$ and $v(\cdot)\in {\cal{V}}[0,T]$. Then, there is a positive constant $C_p$ only depending on  $p$, $T$ and $L$, such that
\begin{eqnarray}\label{fbjia1}
                \mathbb{E}\left[\sup_{t\leq s\leq T}|X^{\gamma_t,u,v}(s)-X^{\gamma_t',u,v}(s)|^p\right]\leq C_p||\gamma_t-\gamma_t'||_0^p;
\end{eqnarray}
\begin{eqnarray}\label{fbjia2}
                \mathbb{E}\left[||X_T^{\gamma_t,u,v}||_0^p\right]\leq C_p(1+||\gamma_t||_0^p);
                \end{eqnarray}
\begin{eqnarray}\label{fbjia3}
                \mathbb{E}\left[||X_{r}^{\gamma_t,u,v}-\gamma_{t}||_0^p\right]\leq C_p{(}1+||\gamma_t||_0^p)(r-t)^{\frac{p}{2}}, \  \ r\in [t,T];
\end{eqnarray}
and
\begin{eqnarray}\label{fbjia4}
                \mathbb{E}\left[\sup_{t\leq s\leq T}|Y^{\gamma_t,u,v}(s)-Y^{\gamma_t',u,v}(s)|^p\right]
                \leq C_p||\gamma_t-\gamma_t'||_0^{p};
                \end{eqnarray}
                \begin{eqnarray}\label{fbjia5}
                \mathbb{E}\left[\sup_{t\leq s\leq T}|Y^{\gamma_t,u,v}(s)|^p\right]
                +\mathbb{E}\left[\left(\int^{T}_{t}|Z^{\gamma_t,u,v}(s)|^2ds\right)^{\frac{p}{2}}\right]\leq C_p(1+||\gamma_t||_0^{p}).
                \end{eqnarray}
\end{lemma}
\par
Formally,  under the assumptions  Hypothesis \ref{hypstate6}, the lower value functional
$\underline{V}(\gamma_t)$ as well as the upper value function $\overline{V}(t, x)$ are  ${\cal{F}}_t$-measurable. 
However, by  the similar proof procedure of Proposition 3.3 in Buckdahn and Li \cite{buck1}, we can prove the following.
\begin{theorem}\label{valuedet} \ \  Suppose the Hypothesis \ref{hypstate6} holds true.
                      Then lower value functional $\underline{V}$ and  upper value functional $\overline{V}$ are deterministic functionals.
\end{theorem}
\par
We now discuss the dynamic programming principle  (DPP) for PSDG (\ref{state1}), (\ref{fbsde1}), (\ref{value1}) and (\ref{value2}).
For this purpose, we define the family of backward semigroups associated with BSDE (\ref{fbsde1}), following the
idea of Peng \cite{peng1}.
\par
Given the initial condition $(t,\gamma_t)\in [0,T)\times{\Lambda}$, a positive number $\delta\leq T-t$, two admissible controls $u(\cdot)\in {\cal{U}}[t,t+\delta]$,  $v(\cdot)\in {\cal{V}}[t,t+\delta]$ and
a real-valued random variable $\eta\in L^2(\Omega,{\cal{F}}_{t+\delta},\mathbb{P};\mathbb{R})$, we put
\begin{eqnarray}\label{gdpp}
                        G^{\gamma_t,u,v}_{s,t+\delta}[\eta]:=\tilde{Y}^{\gamma_t,u,v}(s),\ \
                        \ \ \ \ s\in[t,t+\delta],
\end{eqnarray}
                        where $(\tilde{Y}^{\gamma_t,u,v}(s),\tilde{Z}^{\gamma_t,u,v}(s))_{t\leq s\leq
                        t+\delta}$ is the solution of the following
                        BSDE with the time horizon $t+\delta$:
\begin{eqnarray}\label{bsdegdpp}
\begin{cases}
d\tilde{Y}^{\gamma_t,u,v}(s) =-q(X^{\gamma_t,u,v}_s,\tilde{Y}^{\gamma_t,u,v}(s),\tilde{Z}^{\gamma_t,u,v}(s),u(s),v(s))ds+\tilde{Z}^{\gamma_t,u,v}(s)dW(s),\\
 ~\tilde{Y}^{\gamma_t,u,v}(t+\delta)=\eta,
 \end{cases}
\end{eqnarray}
                  and $X^{\gamma_t,u,v}(\cdot)$ is the solution of PSDE (\ref{state1}).
\begin{theorem}\label{theoremddp} (\cite[Theorem 2.9 ]{zhang})
    Assume Hypothesis \ref{hypstate6}  holds true, the lower value functional
                              $\underline{V}$  and upper value functional $\overline{V}$ obey the following DPPs: for
                              any $(t,\gamma_t)\in [0,T)\times{\Lambda}$ and $0<\delta\leq T-t$,
\begin{eqnarray}\label{ddpG}
                              \underline{V}(\gamma_t)=\mathop{\essinf}\limits_{\beta(\cdot)\in {\cal{B}}_{[t,t+\delta]}}\mathop{\esssup}\limits_{u(\cdot)\in{\mathcal
                              {U}}[t,t+\delta]}G^{\gamma_t,u,\beta(u)}_{t,t+\delta}\left[\underline{V}(X^{\gamma_t,u,\beta(u)}_{t+\delta})\right];
\end{eqnarray}
\begin{eqnarray}\label{ddpG1}
                              \overline{V}(\gamma_t)=\mathop{\esssup}\limits_{\alpha(\cdot)\in {\cal{A}}_{[t,t+\delta]}}\mathop{\essinf}\limits_{v(\cdot)\in{\mathcal
                              {V}}[t,t+\delta]}G^{\gamma_t,\alpha(v),v}_{t,t+\delta}\left[\overline{V}(X^{\gamma_t,\alpha(v),v}_{t+\delta})\right].
\end{eqnarray}
\end{theorem}
\par
From Lemma \ref{lemmaexist} and Theorem \ref{theoremddp}, it follows that the regularity of the value functionals $\underline{V}$ and $\overline{V}$.  
\begin{lemma}\label{lemmavaluev} (\cite[Lemma 2.8 and Proposition 2.13]{zhang})
\ \ Assume that Hypothesis  \ref{hypstate6}     holds, then 
 there is a constant $C>0$ such that, for every  $0\leq t\leq s\leq T$ and $\gamma_t, \eta_{s}\in {\Lambda}$,
\begin{eqnarray}\label{valuelip}
               |\underline{V}(\gamma_t)|\vee|\overline{V}(\gamma_t)|\leq C(1+||\gamma_t||_0).
\end{eqnarray}
\begin{eqnarray}\label{hold}
                 |\underline{V}(\gamma_t)-\underline{V}(\eta_{s})|\vee|\overline{V}(\gamma_t)-\overline{V}(\eta_{s})|\leq
                                    C[||\gamma_{t}-\eta_{s}||_0+(1+||\gamma_t||_0)(s-t)^{\frac{1}{2}}].
\end{eqnarray}
\end{lemma}
We are now in a  position  to give  the existence  result 
for  viscosity solutions. The proof will be given in Appendix B.
\begin{theorem}\label{theoremvexist} \ \
                          Suppose that Hypothesis \ref{hypstate6}  holds. Then  $\underline{V}$ (resp., $\overline{V}$)  is a
                          viscosity solution to equation  (\ref{hjb106}) (resp., (\ref{hjb2})).
\end{theorem}
Hypothesis \ref{hypstate} (i) follows from the continuity of $b(\cdot, u,v)$, $\sigma(\cdot,u,v)$, $q(\cdot,\cdot,\cdot,u,v)$, uniform in $(u,v)\in U\times V$. For Hypothesis \ref{hypstate} (iii)
   we can argue as follows: since $\sigma(\gamma_r,u,v)\sigma^\top(\gamma_r,u,v)$ is a nonnegative, self-adjoint, trace class operator, it is obvious that, for $X,Y\in \mathcal{S}(\mathbb{R}^d)$ with $X\leq Y$,
 $$
 \mbox{tr}[ X \sigma(\gamma_t,u,v)\sigma^\top(\gamma_t,u,v)]\leq \mbox{tr}[ Y \sigma(\gamma_t,u,v)\sigma^\top(\gamma_t,u,v)],
 $$
 and then taking the infimum over $v\in V$ and supremum over $u\in U$ we see that $\mathbf{H}^-$ satisfies Hypothesis \ref{hypstate} (iii). Similarly, taking the supremum over $u\in U$ and  infimum over $v\in V$  we see that $\mathbf{H}^+$ satisfies Hypothesis \ref{hypstate} (iii).
 \par
 To show that Hypothesis \ref{hypstate} (iv) holds observe that, using Hypothesis \ref{hypstate6}, for every $(t,\gamma_t,\eta_t,r)\in [0,T]\times\Lambda\times\Lambda\times  \mathbb{R}$, 
    for any
$\beta>0$, for all $X,Y\in {\cal{S}}(\mathbb{R}^d)$ satisfying
 \begin{eqnarray*}
                              -3\beta\left(\begin{array}{cc}
                                    I&0\\
                                    0&I
                                    \end{array}\right)\leq \left(\begin{array}{cc}
                                    X&0\\
                                    0&Y
                                    \end{array}\right)\leq  3\beta \left(\begin{array}{cc}
                                    I&-I\\
                                    -I&I
                                    \end{array}\right),
\end{eqnarray*}
   we have
  \begin{eqnarray*}
              &&\mathbf{H}^-(\gamma_t,r,\beta(\gamma_t(t)-\eta_t(t)),X)-\mathbf{H}^-(\eta_t,r,\beta(\gamma_t(t)-\eta_t(t)),-Y)\\
              &\leq& \sup_{(u,v)\in U\times V}\bigg{[}( {b}(\gamma_{t},u,v)  -{b}(\eta_{t},u,v),\beta(\gamma_t(t)-\eta_t(t)))_{\mathbb{R}^d}
              +\frac{1}{2}\mbox{tr}{[}X{\sigma}(\gamma_{t},u,v)
                                        {\sigma}^\top(\gamma_{t},u,v){]}\\
                                        &&+\frac{1}{2}\mbox{tr}{[}Y\sigma(\eta_{t},u)\sigma^\top(\eta_{t},u,v){]}
                                        +q{(}\gamma_{t}, r, \sigma^\top(\gamma_{t},u,v)\beta(\gamma_t(t)-\eta_t(t)),u,v{)}\\
                                        &&-
                                 q{(}\eta_{t}, r, \sigma(\eta_{t},u)^\top\beta(\gamma_t(t)-\eta_t(t)),u,v{)}\bigg{]}\\
              &\leq& \beta L||\gamma_t-\eta_t||^2_0 +\frac{3}{2}\beta L^2||\gamma_t-\eta_t||^2_0 +L||\gamma_t-\eta_t||_0+\beta L^2||\gamma_t-\eta_t||^2_0.
\end{eqnarray*}
Taking $\rho(s,r)=(L+3L^2)s$,   Hypothesis \ref{hypstate} (iv) is satisfied for $\mathbf{H}^-$. Similarly, $\mathbf{H}^+$ satisfies Hypothesis \ref{hypstate} (iv). Hypothesis \ref{hypstate} (v)
follows from (\ref{211224}) and (\ref{2112243}).    
By    \cite[Proposition 11.2.13]{zhang1},  without loss of generality we may assume that
there exists a constant $\nu\geq0$
                 such that, for every $(t,\gamma_t,p,X)
                 \in [0,T]\times\Lambda\times   \mathbb{R}^d\times {\cal{S}}(\mathbb{R}^d)$,
      \begin{eqnarray*}
              ({\mathbf{H}}^-(\gamma_t,r,p,X)-{\mathbf{H}}^-(\gamma_t,s,p,X))\vee  ({\mathbf{H}}^+(\gamma_t,r,p,X)-{\mathbf{H}}^+(\gamma_t,s,p,X))\geq
                 \nu(s-r) \ \  \mbox{when} \  r\leq s.
\end{eqnarray*}
                   Then    Theorems    \ref{theoremvexist} and \ref{theoremhjbm} lead to the result (given below) that the viscosity solution to   PHJBIE given in  (\ref{hjb106}) (resp., (\ref{hjb2}))
                      corresponds to the lower value functional  $\underline{V}$ (resp.,  upper value functional $\overline{V}$) of our PSDG given in (\ref{state1}), (\ref{fbsde1}) and (\ref{value1}) (resp., (\ref{value2})).
\begin{theorem}\label{theorem52}\ \
                 Let Hypothesis \ref{hypstate6}  hold. Then $\underline{V}$ (resp., $\overline{V}$) defined by (\ref{value1}) (resp., (\ref{value2})) is the unique viscosity
                          solution to (\ref{hjb106}) (resp., (\ref{hjb2})) in the class of functionals satisfying (\ref{w}) and (\ref{w1}).
\end{theorem}

\begin{remark}\label{remarkv12241}
If the Isaacs' condition holds, that is, if for all $(t,\gamma_t,r,p,X)
                 \in [0,T]\times\Lambda\times \mathbb{R}\times  \mathbb{R}^d\times {\cal{S}}(\mathbb{R}^d)$,
$$
{\mathbf{H}}^-(\gamma_t,r,p,X)={\mathbf{H}}^+(\gamma_t,r,p,X),
$$
then (\ref{hjb106}) and (\ref{hjb2}) coincide, and from the uniqueness  of viscosity solutions it
follows that the lower value functional $\underline{V}$ equals the upper value functional $\overline{V}$
which means that the associated stochastic differential game has a value.
\end{remark}

\appendix

   \section{  Consistency and stability for viscosity solutions.}

\setcounter{equation}{0}
\renewcommand{\theequation}{A.\arabic{equation}}
{\bf  Proof of Lemma \ref{theorem3.20407}}. \ \   Let  $\varphi\in \mathcal{A}^+(\gamma_t,v)$ with $(t,{\gamma}_{t})\in [\hat{t},T)\times B_r(\hat{\gamma}_{\hat{t}})$.
 For  every $\alpha\in \mathbb{R}^{d}$ and $\beta\in \mathbb{R}^{d\times n}$,
 let  $$
           X(s)=\gamma_t(t)+\int^{s}_{t}\alpha dl+\int^{s}_{t}\beta dW(l),\ \ s\in [t,T],
 $$
 and $X(s)=\gamma_t(s)$, $s\in [0,t)$. Then $X(\cdot)$ is a continuous semi-martingale on $[t,T]$.
  Applying functional It\^o formula (\ref{statesop}) to $\varphi$ and noticing that $(v-\varphi)(\gamma_t)=0$, we have, for every $0<\delta\leq T-t$,
 \begin{eqnarray}\label{1224}
                            0&\leq&{\mathbb{E}}(\varphi-v)(X_{t+\delta})\nonumber\\
                            &=&{\mathbb{E}}\int^{t+\delta}_{t}[\partial_t(\varphi-v)(X_l)
                 +(\partial_x(\varphi-v)(X_l),\alpha)_{\mathbb{R}^d}]dl+\frac{1}{2}{\mathbb{E}}\int^{t+\delta}_{t}\mbox{tr}[(\partial_{xx}(\varphi-v)(X_l))\beta\beta^\top]dl\nonumber\\
                 &=&
                           {\mathbb{E}}\int^{t+\delta}_{t}\widetilde{{\mathcal{H}}}(X_l)
                              dl,
\end{eqnarray}
where \begin{eqnarray*}
\widetilde{{\mathcal{H}}}(\eta_s)=\partial_t(\varphi-v)(\eta_s)
                 +(\partial_x(\varphi-v)(\eta_s),\alpha)_{\mathbb{R}^d}
                 +\frac{1}{2}\mbox{tr}[(\partial_{xx}(\varphi-v)(\eta_s))\beta\beta^\top], \ \ \ (s,\eta_s)\in [0,T]\times \Lambda.
                 \end{eqnarray*}
                 Letting $\delta\rightarrow0$ in (\ref{1224}),
\begin{eqnarray}\label{040912}
\widetilde{{\mathcal{H}}}(\gamma_t)\geq0.
\end{eqnarray}
Let $\beta=\mathbf{0}$, by the arbitrariness of  $\alpha$,
$$
\partial_t\varphi(\gamma_t)\geq\partial_tv(\gamma_t), \ \ \partial_x\varphi(\gamma_t)=\partial_xv(\gamma_t). 
$$
Let  $\alpha=\mathbf{0}$,
$$
\frac{1}{2}\mbox{tr}[(\partial_{xx}(\varphi-v)(\gamma_t))\beta\beta^\top]+\partial_t(\varphi-v)(\gamma_t)\geq0.
$$
 By the arbitrariness of  $\beta$,
 $$
 \partial_{xx}\varphi(\gamma_t)\geq \partial_{xx}v(\gamma_t).
 $$
Note that $\mathcal{L}v(\gamma_t)\geq0$, by Hypothesis \ref{hypstate} (iii),
we have
$
\mathcal{L}\varphi(\gamma_t)\geq\mathcal{L}v(\gamma_t)\geq0
$.
 Thus,
   $v$ is a viscosity subsolution of PPDE (\ref{hjb1}). 
    \ \ $\Box$
    \par
{\bf  Proof of Theorem \ref{theorem3.2}}. \ \  We prove the subsolution property only.  Assume $v$ is a viscosity subsolution. It is clear that $v(\gamma_T)\leq\phi(\gamma_T)$ for all $\gamma_T\in \Lambda_T$. For any $(t,\gamma_t)\in [0,T)\times \Lambda$, since $v\in C_p^{1,2}({\Lambda})$, by definition of viscosity subsolutions we see that $\mathcal{L}v(\gamma_t)\geq 0$.

On the other hand, assume $v$ is a  classical subsolution on $\Lambda$.    Let  $\varphi\in \mathcal{A}^+(\gamma_t,v)$ with $t\in [0,T)$. By the same proof procedure of Lemma  \ref{theorem3.20407}, we have
$
\mathcal{L}\varphi(\gamma_t)\geq\mathcal{L}v(\gamma_t)\geq0
$ and
   $v$ is a viscosity subsolution of PPDE (\ref{hjb1}). 
    \ \ $\Box$

\par
{\bf  Proof of Theorem \ref{theoremstability}}. \ \ Without loss of generality, we shall only prove the viscosity subsolution property.
First, since $v^{\varepsilon}$ is a viscosity subsolution of  PPDE (\ref{hjb1}) with generator $\mathbf{F}^{\varepsilon}$, we have
 $$
                            v^{\varepsilon}(\gamma_T)\leq \phi^{\varepsilon}(\gamma_T),\ \ \gamma_T\in \Lambda_T.
$$
Letting $\varepsilon\rightarrow0$,
 $$
                            v(\gamma_T)\leq \phi(\gamma_T),\ \ \gamma_T\in \Lambda_T.
$$
Next,      let   $\varphi\in \mathcal{A}^+(\hat{\gamma}_{\hat{t}}, v)$ with
  $(\hat{t},\hat{\gamma}_{\hat{t}})\in [0,T)\times\Lambda$. By (\ref{sss}), 
   there exists a constant $\Delta>0$ such that for all $\varepsilon\in (0,\Delta)$,
  $$\sup_{(t,\gamma_t)\in [\hat{t},T]\times\Lambda^{\hat{t}}}(v^\varepsilon(\gamma_t)-\varphi(\gamma_t))\leq 1.$$
 Denote $\varphi_{1}(\gamma_t):=\varphi(\gamma_t)+\overline{\Upsilon}_0(\gamma_t,\hat{\gamma}_{{\hat{t}}})$
 for all
 $(t,\gamma_t)\in [0,T]\times\Lambda$. By Lemma \ref{theoremS}, we have  $\varphi_{1}\in C^{1,2}_{p}(\Lambda^{\hat{t}})$.
 For every  $\varepsilon\in (0,\Delta)$, it is clear that $v^{\varepsilon}-{{\varphi_{1}}}$ is an  upper semicontinuous functional
  and bounded from above on $\Lambda^{\hat{t}}$.
  Define a sequence of positive numbers $\{\delta_i\}_{i\geq0}$  by 
        $\delta_i=\frac{1}{2^i}$ for all $i\geq0$. 
        Since $\overline{\Upsilon}_0(\cdot,\cdot)$ is a gauge-type function on compete metric space $(\Lambda^{\hat{t}},d_{\infty})$, from Lemma \ref{theoremleft} it follows that,
 for every  $(t_0,\gamma^0_{t_0})\in [\hat{t},T]\times \Lambda^{\hat{t}}$ satisfying
\begin{eqnarray}\label{0513a}
(v^{\varepsilon}-{{\varphi_{1}}})(\gamma^0_{t_0})\geq \sup_{(s,\gamma_s)\in [\hat{t},T]\times \Lambda^{\hat{t}}}(v^{\varepsilon}-{{\varphi_{1}}})(\gamma_s)-\varepsilon,\
\    \mbox{and} \ \ (v^{\varepsilon}-{{\varphi_{1}}})(\gamma^0_{t_0})\geq (v^{\varepsilon}-{{\varphi_{1}}})(\hat{\gamma}_{\hat{t}}),
\end{eqnarray}
  there exist $(t_{\varepsilon},{\gamma}^{\varepsilon}_{t_{\varepsilon}})\in [\hat{t},T]\times \Lambda^{\hat{t}}$ and a sequence $\{(t_i,\gamma^i_{t_i})\}_{i\geq1}\subset [\hat{t},T]\times \Lambda^{\hat{t}}$ such that
  \begin{description}
        \item{(i)} $\overline{\Upsilon}_0(\gamma^0_{t_0},{\gamma}^{\varepsilon}_{t_{\varepsilon}})\leq {\varepsilon}$,  $\overline{\Upsilon}_0(\gamma^i_{t_i},{\gamma}^{\varepsilon}_{t_{\varepsilon}})\leq \frac{\varepsilon}{2^i}$ and $t_i\uparrow t_{\varepsilon}$ as $i\rightarrow\infty$,
        \item{(ii)}  $(v^{\varepsilon}-{{\varphi_{1}}})({\gamma}^{\varepsilon}_{t_{\varepsilon}})-\sum_{i=0}^{\infty}\frac{1}{2^i}\overline{\Upsilon}_0(\gamma^i_{t_i},{\gamma}^{\varepsilon}_{t_{\varepsilon}})\geq (v^{\varepsilon}-{{\varphi_{1}}})(\gamma^0_{t_0})$, and
        \item{(iii)}  $(v^{\varepsilon}-{{\varphi_{1}}})(\gamma_s)-\sum_{i=0}^{\infty}\frac{1}{2^i}\overline{\Upsilon}_0(\gamma^i_{t_i},\gamma_s)
            <(v^{\varepsilon}-{{\varphi_{1}}})({\gamma}^{\varepsilon}_{t_{\varepsilon}})-\sum_{i=0}^{\infty}\frac{1}{2^i}\overline{\Upsilon}_0(\gamma^i_{t_i},{\gamma}^{\varepsilon}_{t_{\varepsilon}})$ for all $(s,\gamma_s)\in [t_{\varepsilon},T]\times \Lambda^{t_{\varepsilon}}\setminus \{(t_{\varepsilon},{\gamma}^{\varepsilon}_{t_{\varepsilon}})\}$.

        \end{description}
  %
%
%
%
We claim that
\begin{eqnarray}\label{gamma}
d_\infty({\gamma}^{\varepsilon}_{{t}_{\varepsilon}},\hat{\gamma}_{\hat{t}})\rightarrow0  \ \ \mbox{as} \ \ \varepsilon\rightarrow0.
\end{eqnarray}
 Indeed, if not,  by (\ref{up}) and the definition of $d_{\infty}$, we can assume that there exists a constant  $\nu_0>0$
 such
                    that
$$
  \overline{\Upsilon}_0({\gamma}^{\varepsilon}_{{t}_{\varepsilon}},\hat{\gamma}_{{\hat{t}}})
  \geq\nu_0.
$$
 Thus, by (\ref{0513a}) and  the property (ii) of $(t_{\varepsilon},{\gamma}^{\varepsilon}_{t_{\varepsilon}})$,  we obtain that
\begin{eqnarray*}
   &&0=(v- {{\varphi}})(\hat{\gamma}_{\hat{t}})= \lim_{\varepsilon\rightarrow0}(v^\varepsilon-{{\varphi_{1}}})(\hat{\gamma}_{\hat{t}})
   \leq {\limsup_{\varepsilon\rightarrow0}}\bigg{[}(v^{\varepsilon}-{{\varphi_{1}}})({\gamma}^{\varepsilon}_{t_{\varepsilon}})-\sum_{i=0}^{\infty}\frac{1}{2^i}\overline{\Upsilon}_0(\gamma^i_{t_i},{\gamma}^{\varepsilon}_{t_{\varepsilon}})\bigg{]}\\
   &=&{\limsup_{\varepsilon\rightarrow0}}\bigg{[}(v^\varepsilon-{{\varphi}})({\gamma}^{\varepsilon}_{{t}_{\varepsilon}})
  -\overline{\Upsilon}_0({\gamma}^{\varepsilon}_{{t}_{\varepsilon}},\hat{\gamma}_{{\hat{t}}})-
  \sum_{i=0}^{\infty}\frac{1}{2^i}\overline{\Upsilon}_0(\gamma^i_{t_i},{\gamma}^{\varepsilon}_{t_{\varepsilon}})\bigg{]}\\
   &\leq&{\limsup_{\varepsilon\rightarrow0}}\bigg{[}(v-{{\varphi}})({\gamma}^{\varepsilon}_{{t}_{\varepsilon}})+(v^\varepsilon-v)({\gamma}^{\varepsilon}_{{t}_{\varepsilon}})
     -\sum_{i=0}^{\infty}\frac{1}{2^i}\overline{\Upsilon}_0(\gamma^i_{t_i},{\gamma}^{\varepsilon}_{t_{\varepsilon}})\bigg{]}-\nu_0\leq (v- {{\varphi}})(\hat{\gamma}_{\hat{t}})-\nu_0=-\nu_0,
\end{eqnarray*}
 contradicting $\nu_0>0$.  We notice that, by (\ref{s0}) and the property (i) of $(t_{\varepsilon},{\gamma}^{\varepsilon}_{t_{\varepsilon}})$, exists a generic constant $C>0$ such that
 \begin{eqnarray*}
  2\sum_{i=0}^{\infty}\frac{1}{2^i}({t_{\varepsilon}}-{t}_{i})
  \leq2\sum_{i=0}^{\infty}\frac{1}{2^i}\bigg{(}\frac{\varepsilon}{2^i}\bigg{)}^{\frac{1}{2}}\leq C\varepsilon^{\frac{1}{2}};
    \end{eqnarray*}
    \begin{eqnarray*}
    |{\partial_x{\Upsilon}}({\gamma}^{\varepsilon}_{{t}_{\varepsilon}}-\hat{\gamma}_{{\hat{t}},{t}_{\varepsilon}})|\leq C|\hat{\gamma}_{{\hat{t}}}(\hat{t})-{\gamma}^{\varepsilon}_{{t}_{\varepsilon}}({t}_{\varepsilon})|^5;\ \
   |{\partial_{xx}{\Upsilon}}({\gamma}^{\varepsilon}_{{t}_{\varepsilon}}-\hat{\gamma}_{{\hat{t}},{t}_{\varepsilon}})|\leq C|\hat{\gamma}_{{\hat{t}}}(\hat{t})-{\gamma}^{\varepsilon}_{{t}_{\varepsilon}}({t}_{\varepsilon})|^4;
    \end{eqnarray*}
    \begin{eqnarray*}
  \left|\partial_x\left[\sum_{i=0}^{\infty}\frac{1}{2^i}
                      \Upsilon({\gamma}^{\varepsilon}_{t_{\varepsilon}}-\gamma^i_{t_i,t_{\varepsilon}})\right]
                      \right|
                      \leq18\sum_{i=0}^{\infty}\frac{1}{2^i}|\gamma^i_{t_i}({t}_{i})-{\gamma}^{\varepsilon}_{t_{\varepsilon}}(t_{\varepsilon})|^5
                     \leq18\sum_{i=0}^{\infty}\frac{1}{2^i}\bigg{(}\frac{\varepsilon}{2^i}\bigg{)}^{\frac{5}{6}}
                      \leq C{\varepsilon}^{\frac{5}{6}};
                        \end{eqnarray*}
    and
     \begin{eqnarray*}
 \left|\partial_{xx}\left[\sum_{i=0}^{\infty}\frac{1}{2^i}
                      \Upsilon({\gamma}^{\varepsilon}_{t_{\varepsilon}}-\gamma^i_{t_i,t_{\varepsilon}})\right]
                      \right|
                      \leq306\sum_{i=0}^{\infty}\frac{1}{2^i}|\gamma^i_{t_i}({t}_{i})-{\gamma}^{\varepsilon}_{t_{\varepsilon}}(t_{\varepsilon})|^4
                     \leq306\sum_{i=0}^{\infty}\frac{1}{2^i}\bigg{(}\frac{\varepsilon}{2^i}\bigg{)}^{\frac{2}{3}}
                      \leq C{\varepsilon}^{\frac{2}{3}}.
                        \end{eqnarray*}
 Then for any $\varrho>0$, by (\ref{sss}), (\ref{gamma}) and Hypothesis \ref{hypstate} (i), there exists $\varepsilon>0$ small enough such that
$$
            \hat{t}\leq {t}_{\varepsilon}< T,  \        \  
             2|{t}_{\varepsilon}-\hat{t}|+2\sum_{i=0}^{\infty}\frac{1}{2^i}({t_{\varepsilon}}-{t}_{i})\leq \frac{\varrho}{3}, $$
             and
             $$
|\partial_t{\varphi}({\gamma}^{\varepsilon}_{{t}_{\varepsilon}})-\partial_t{\varphi}(\hat{\gamma}_{\hat{t}})|\leq \frac{\varrho}{3}, \ |I|\leq \frac{\varrho}{3}, 
$$
where
\begin{eqnarray*}
I:={\mathbf{F}}^{\varepsilon}({\gamma}^{\varepsilon}_{{t}_{\varepsilon}}, v^{\varepsilon}({\gamma}^{\varepsilon}_{{t}_{\varepsilon}}),
                           \partial_x{\varphi_{2}}({\gamma}^{\varepsilon}_{{t}_{\varepsilon}}),\partial_{xx}{\varphi_{2}}({\gamma}^{\varepsilon}_{{t}_{\varepsilon}}))
                       -{\mathbf{F}}(\hat{\gamma}_{\hat{t}},{\varphi}(\hat{\gamma}_{\hat{t}}),\partial_x{\varphi}(\hat{\gamma}_{\hat{t}}),\partial_{xx}{\varphi}(\hat{\gamma}_{\hat{t}})),
\end{eqnarray*}
and
\begin{eqnarray*} {\varphi_{2}}({\gamma}^{\varepsilon}_{{t}_{\varepsilon}})={\varphi_{1}}({\gamma}^{\varepsilon}_{{t}_{\varepsilon}})+\sum_{i=0}^{\infty}\frac{1}{2^i}\overline{\Upsilon}_0(\gamma^i_{t_i},{\gamma}^{\varepsilon}_{t_{\varepsilon}}).
\end{eqnarray*}
 Since $v^{\varepsilon}$ is a viscosity subsolution of PPDE (\ref{hjb1}) with generators $\mathbf{F}^{\varepsilon}$, we have
$$
                          \partial_t{\varphi_{2}}({\gamma}^{\varepsilon}_{{t}_{\varepsilon}})
                           +{\mathbf{F}}^{\varepsilon}({\gamma}^{\varepsilon}_{{t}_{\varepsilon}}, v^{\varepsilon}({\gamma}^{\varepsilon}_{{t}_{\varepsilon}}),
                           \partial_x{\varphi_{2}}({\gamma}^{\varepsilon}_{{t}_{\varepsilon}}),\partial_{xx}{\varphi_{2}}({\gamma}^{\varepsilon}_{{t}_{\varepsilon}}))\geq0.
$$
Thus
\begin{eqnarray*}
                       0\leq  \partial_t{\varphi}({\gamma}^{\varepsilon}_{{t}_{\varepsilon}})
                       +2({t}_{\varepsilon}-\hat{t})+2\sum_{i=0}^{\infty}({t_{\varepsilon}}-{t}_{i})
                       +{\mathbf{F}}(\hat{\gamma}_{\hat{t}},{\varphi}(\hat{\gamma}_{\hat{t}}),\partial_x{\varphi}(\hat{\gamma}_{\hat{t}}),
                       \partial_{xx}{\varphi}(\hat{\gamma}_{\hat{t}}))+I\leq \mathcal{L}v(\hat{\gamma}_{\hat{t}})+\varrho.
\end{eqnarray*}
Letting $\varrho\downarrow 0$, we show that $\mathcal{L}v(\hat{\gamma}_{\hat{t}})\geq0$ and
 $v$ is a viscosity subsolution of  PPDE (\ref{hjb1}) with generator $\mathbf{F}$. 
\ \ $\Box$

 \section{  Existence  for viscosity solutions to PHJBIEs.}

\setcounter{equation}{0}
\renewcommand{\theequation}{B.\arabic{equation}}

{\bf  Proof of Theorem \ref{theoremvexist}}. \ \ We shall only prove that $\underline{V}$ satisfies viscosity property of equation  (\ref{hjb106}). The other statements can be proved similarly.
  First, let us show that $\underline{V}$ is a viscosity subsolution of (\ref{hjb106}).
 We let  $\varphi\in \mathcal{A}^+(\hat{\gamma}_{\hat{t}},\underline{V})$
                  with
                   $(\hat{t},\hat{\gamma}_{\hat{t}})\in [0,T)\times \Lambda$.  We need to prove that
\begin{eqnarray}\label{628jia0}
       \partial_t{\varphi}(\hat{\gamma}_{\hat{t}})
                           +{\mathbf{H}}^-(\hat{\gamma}_{\hat{t}},{\varphi}(\hat{\gamma}_{\hat{t}}),\partial_x{\varphi}(\hat{\gamma}_{\hat{t}}),\partial_{xx}{\varphi}(\hat{\gamma}_{\hat{t}}))\geq0.
\end{eqnarray}
     Thus, we assume the converse result that $\theta>0$ exists  such that
       $$
       \partial_t{\varphi}(\hat{\gamma}_{\hat{t}})
                           +{\mathbf{H}}^-(\hat{\gamma}_{\hat{t}},{\varphi}(\hat{\gamma}_{\hat{t}}),\partial_x{\varphi}(\hat{\gamma}_{\hat{t}}),\partial_{xx}{\varphi}(\hat{\gamma}_{\hat{t}}))\leq-2\theta<0.
 $$
 Then by the following Lemma B.1. we can find a measurable function $\pi:U\rightarrow V$ such that
\begin{eqnarray}\label{628jia}
 \partial_t{\varphi}(\hat{\gamma}_{\hat{t}})
                           +{\mathbf{H}}(\hat{\gamma}_{\hat{t}},{\varphi}(\hat{\gamma}_{\hat{t}}),\partial_x{\varphi}(\hat{\gamma}_{\hat{t}}),\partial_{xx}{\varphi}(\hat{\gamma}_{\hat{t}}),u,\pi(u))\leq-\theta,\ \ \mbox{for all}\ u\in U.
\end{eqnarray}
 On the other hand, for  $0< \delta\leq T-\hat{t}$, we have $\hat{t}< \hat{t}+\delta \leq T$, then by the DPP (Theorem \ref{theoremddp}), we obtain the following result:
 \begin{eqnarray}\label{4.9}
                           && 0=\underline{V}(\hat{\gamma}_{\hat{t}})-{{\varphi}} (\hat{\gamma}_{\hat{t}})
                           =\mathop{\essinf}\limits_{\beta(\cdot)\in {\cal{B}}_{[\hat{t},\hat{t}+\delta]}}\mathop{\esssup}\limits_{u(\cdot)\in{\mathcal
                              {U}}[\hat{t},\hat{t}+\delta]}G^{\hat{\gamma}_{\hat{t}},u,\beta(u)}_{\hat{t},\hat{t}+\delta}[\underline{V}(X^{\hat{\gamma}_{\hat{t}},u,\beta(u)}_{\hat{t}+\delta})]
                           -{{\varphi}} (\hat{\gamma}_{\hat{t}}).
\end{eqnarray}
In particular,
\begin{eqnarray}\label{4.9jia628}
                           && 0=\underline{V}(\hat{\gamma}_{\hat{t}})-{{\varphi}} (\hat{\gamma}_{\hat{t}})
                           \leq\mathop{\esssup}\limits_{u(\cdot)\in{\mathcal
                              {U}}[\hat{t},\hat{t}+\delta]}G^{\hat{\gamma}_{\hat{t}},u,\pi(u)}_{\hat{t},\hat{t}+\delta}[\underline{V}(X^{\hat{\gamma}_{\hat{t}},u,\pi(u)}_{\hat{t}+\delta})]
                           -{{\varphi}} (\hat{\gamma}_{\hat{t}}).
\end{eqnarray}
Here, by putting $\pi(u)(s,\omega)=\pi(u(s,\omega))$, $(s,\omega)\in [\hat{t},T]\times \Omega$, we identify $\pi$ as an element of ${\mathcal{B}}_{[\hat{t},\hat{t}+\delta]}$.
Then, for any $\varepsilon>0$ and $0<\delta\leq T-\hat{t}$,  we can  find a control  ${u}^{\varepsilon}(\cdot)
\equiv u^{{\varepsilon},\delta}(\cdot)
\in {\cal{U}}[\hat{t},\hat{t}+\delta]$ such
   that 
   the following result holds:
\begin{eqnarray}\label{4.10}
    -{\varepsilon}\delta
    \leq G^{\hat{\gamma}_{\hat{t}},{u}^{\varepsilon},\pi({u}^{\varepsilon})}_{{\hat{t}},\hat{t}+\delta}[\underline{V}(X^{\hat{\gamma}_{\hat{t}},{u}^{\varepsilon},\pi({u}^{\varepsilon})}_{{\hat{t}}+\delta})]-{{\varphi}} (\hat{\gamma}_{\hat{t}}).
\end{eqnarray}
 We note that
                     $G^{\hat{\gamma}_{\hat{t}},{u}^{\varepsilon},\pi({u}^{\varepsilon})}_{s,\hat{t}+\delta}[\underline{V}(X^{\hat{\gamma}_{\hat{t}},{u}^{\varepsilon},\pi({u}^{\varepsilon})}_{{\hat{t}}+\delta})]$
                     is defined in terms of the solution of the
                     BSDE:
 \begin{eqnarray}\label{bsde4.10}
\begin{cases}
dY^{\hat{\gamma}_{\hat{t}},{u}^{\varepsilon},\pi({u}^{\varepsilon})}(s) =
               -q(X^{\hat{\gamma}_{\hat{t}},{u}^{\varepsilon},\pi({u}^{\varepsilon})}_s,Y^{\hat{\gamma}_{\hat{t}},{u}^{\varepsilon},\pi({u}^{\varepsilon})}(s),Z^{\hat{\gamma}_{\hat{t}},
               {u}^{\varepsilon}(s),\pi({u}^{\varepsilon})}(s),{u}^{\varepsilon}(s),\pi({u}^{\varepsilon})(s))ds\\
               ~~~~~~~~~~~~~~~~~~~~~~~+Z^{\hat{\gamma}_{\hat{t}},{u}^{\varepsilon},\pi({u}^{\varepsilon})}(s)dW(s),\ \  s\in[\hat{t},\hat{t}+\delta], \\
 ~Y^{\hat{\gamma}_{\hat{t}},{u}^{\varepsilon},\pi({u}^{\varepsilon})}(\hat{t}+\delta)=\underline{V}(X^{\hat{\gamma}_{\hat{t}},
 {u}^{\varepsilon},\pi({u}^{\varepsilon})}_{{\hat{t}}+\delta}),
 \end{cases}
\end{eqnarray}
                     by the following formula:
$$
                        G^{\hat{\gamma}_{\hat{t}},{u}^{\varepsilon},\pi({u}^{\varepsilon})}_{s,\hat{t}+\delta}[\underline{V}(X^{\hat{\gamma}_{\hat{t}},{u}^{\varepsilon},\pi({u}^{\varepsilon})}_{{\hat{t}}+\delta})]
                        =Y^{\hat{\gamma}_{\hat{t}},{u}^{\varepsilon},\pi({u}^{\varepsilon})}(s),\
                        \ \ s\in[\hat{t},\hat{t}+\delta].
$$
                         Applying functional It\^{o} formula (\ref{statesop}) to
                         ${\varphi}(X^{\hat{\gamma}_{\hat{t}},{u}^{\varepsilon},\pi({u}^{\varepsilon})}_s)$,   we get that
\begin{eqnarray}\label{bsde4.21}
                             &&{\varphi}(X^{\hat{\gamma}_{\hat{t}},{u}^{\varepsilon},\pi({u}^{\varepsilon})}_s)\nonumber\\
                             &=& {\varphi}(\hat{\gamma}_{\hat{t}})+\int^{s}_{{\hat{t}}} (\tilde{\cal{L}}{\varphi})(X^{\hat{\gamma}_{\hat{t}},{u}^{\varepsilon},\pi({u}^{\varepsilon})}_l,
                                      u^{{\varepsilon}}(l),\pi(u^{{\varepsilon}})(l))dl -\int^{s}_{{\hat{t}}}q(X^{\hat{\gamma}_{\hat{t}},{u}^{\varepsilon},\pi({u}^{\varepsilon})}_l,{\varphi}(X^{\hat{\gamma}_{\hat{t}},
                                      {u}^{\varepsilon},\pi({u}^{\varepsilon})}_l),\nonumber\\
                                      && \sigma^\top(X^{\hat{\gamma}_{\hat{t}},{u}^{\varepsilon},\pi({u}^{\varepsilon})}_l,u^{{\varepsilon}}(l),\pi(u^{{\varepsilon}})(l))
                                      \partial_x{\varphi}(X^{\hat{\gamma}_{\hat{t}},{u}^{\varepsilon},\pi({u}^{\varepsilon})}_l), u^{\varepsilon}(l),\pi(u^{\varepsilon})(l))dl\nonumber\\
                             &&
                             +\int^{s}_{{\hat{t}}}
                             \sigma^\top(X^{\hat{\gamma}_{\hat{t}},{u}^{\varepsilon},\pi({u}^{\varepsilon})}_l,u^{\varepsilon}(l),\pi(u^{\varepsilon})(l))
                             \partial_x{\varphi}(X^{\hat{\gamma}_{\hat{t}},{u}^{\varepsilon},\pi({u}^{\varepsilon})}_l)dW(l),\ \
\end{eqnarray}
where
\begin{eqnarray*}
                       (\tilde{\cal{L}}{\varphi})(\gamma_t,u,v)
                       &=&\partial_t{\varphi}(\gamma_t)+
                                         (\partial_x {\varphi}(\gamma_t),b(\gamma_t,u,v))_{\mathbb{R}^d}+\frac{1}{2}\mbox{tr}[\partial_{xx}{\varphi}(\gamma_t)\sigma(\gamma_t,u,v)\sigma^{\top}(\gamma_t,u,v)]\\
                       &&+
                                         q(\gamma_t,{\varphi}(\gamma_t),\sigma^{\top}(\gamma_t,u,v){\partial_x{\varphi}(\gamma_t)},u,v),
                                         ~~~~ \ (t, \gamma_t,u,v)\in [0,T]\times {\Lambda}\times U\times V.
\end{eqnarray*}
Set
\begin{eqnarray*}
                          &&Y^{2}(s):=
                          {\varphi}(X_s^{\hat{\gamma}_{\hat{t}},{u}^{\varepsilon},\pi({u}^{\varepsilon})})-Y^{\hat{\gamma}_{\hat{t}},u^{\varepsilon},\pi({u}^{\varepsilon})}(s),  \ \ s\in[\hat{t},\hat{t}+\delta],\\
                            &&Z^{2}(s)  :=
                             \sigma^{\top}(X^{\hat{\gamma}_{\hat{t}},{u}^{\varepsilon},\pi({u}^{\varepsilon})}_s,u^{\varepsilon}(s),\pi(u^{\varepsilon})(s))
                             \partial_x{\varphi}(X^{\hat{\gamma}_{\hat{t}},{u}^{\varepsilon},\pi({u}^{\varepsilon})}_s)-Z^{\hat{\gamma}_{\hat{t}},u^{\varepsilon},\pi({u}^{\varepsilon})}(s), \ \ s\in[\hat{t},\hat{t}+\delta].
\end{eqnarray*}
Comparing (\ref{bsde4.10}) and (\ref{bsde4.21}), we have, $\mathbb{P}$-a.s.,
\begin{eqnarray*}
                      dY^{2}(s)
                      &=&[(\tilde{\cal{L}}{\varphi})(X^{\hat{\gamma}_{\hat{t}},{u}^{\varepsilon},\pi({u}^{\varepsilon})}_s,u^{\varepsilon}(s),\pi(u^{\varepsilon})(s))
                       -q(X^{\hat{\gamma}_{\hat{t}},{u}^{\varepsilon},\pi({u}^{\varepsilon})}_s,{\varphi}(X^{\hat{\gamma}_{\hat{t}},{u}^{\varepsilon},\pi({u}^{\varepsilon})}_s),
                             \\
                             &&
                             \sigma^{\top}(X^{\hat{\gamma}_{\hat{t}},{u}^{\varepsilon},\pi({u}^{\varepsilon})}_s,{u}^{\varepsilon}(s),\pi({u}^{\varepsilon})(s))
                             \partial_x{\varphi}(X^{\hat{\gamma}_{\hat{t}},{u}^{\varepsilon},\pi({u}^{\varepsilon})}_s), {u}^{\varepsilon}(s),\pi({u}^{\varepsilon})(s))\\
                             &&+q(X^{\hat{\gamma}_{\hat{t}},{u}^{\varepsilon},\pi({u}^{\varepsilon})}_s,Y^{\hat{\gamma}_{\hat{t}},{u}^{\varepsilon}(s),\pi({u}^{\varepsilon})}(s),
                                  Z^{\hat{\gamma}_{\hat{t}},{u}^{\varepsilon},\pi({u}^{\varepsilon})}(s),{u}^{\varepsilon}(s),\pi({u}^{\varepsilon})(s))]ds
                                +Z^{2}(s)dW(s)\\
                                &=&[(\tilde{\cal{L}}{\varphi})(X^{\hat{\gamma}_{\hat{t}},{u}^{\varepsilon},\pi({u}^{\varepsilon})}_s,u^{\varepsilon}(s),\pi(u^{\varepsilon})(s))
                                -A(s)Y^{2}(s)-(\bar{A}(s),Z^{2}(s))_{\mathbb{R}^n}]ds
                                +Z^{2}(s)dW(s),
\end{eqnarray*}
where $|A|\vee|\bar{A}|\leq L$. 
By Proposition 2.2 in \cite{el1}, we have
\begin{eqnarray}\label{4.14}
                      Y^{2
                      }(\hat{t})
                      =\mathbb{E}\bigg{[}Y^{2}(t+\delta)\Gamma^{\hat{t}}(\hat{t}+\delta)
                      -
                      \int^{\hat{t}+\delta}_{{\hat{t}}}\Gamma^{\hat{t}}(l)(\tilde{\cal{L}}{\varphi})
                       (X^{\hat{\gamma}_{\hat{t}},{u}^{\varepsilon},\pi({u}^{\varepsilon})}_l,u^{\varepsilon}(l),\pi(u^{\varepsilon})(l))dl\bigg{|}
                      {\cal{F}}_{\hat{t}}\bigg{]},\ \ \
\end{eqnarray}
where $\Gamma^{\hat{t}}(\cdot)$ solves the linear SDE
$$
               d\Gamma^{\hat{t}}(s)=\Gamma^{\hat{t}}(s)(A(s)ds+\bar{A}(s)dW(s)),\ s\in [{\hat{t}},{\hat{t}}+\delta];\ \ \ \Gamma^{\hat{t}}({\hat{t}})=1.
$$
It is clear that  $\Gamma^{\hat{t}}\geq 0$. Combining (\ref{4.10}) and (\ref{4.14}), we obtain
\begin{eqnarray}\label{4.15}
-\varepsilon
    &\leq& \frac{1}{\delta}\mathbb{E}\bigg{[}-Y^{2
    }(\hat{t}+\delta)\Gamma^{\hat{t}}(\hat{t}+\delta)+
    \int^{\hat{t}+\delta}_{{\hat{t}}}\Gamma^{\hat{t}}(l)(\tilde{\cal{L}}{\varphi})(X^{\hat{\gamma}_{\hat{t}},{u}^{\varepsilon},\pi({u}^{\varepsilon})}_l,u^{\varepsilon}(l),\pi(u^{\varepsilon})(l))dl\bigg{]}\nonumber\\
    &=&-\frac{1}{\delta}\mathbb{E}\bigg{[}Y^{2
    }(\hat{t}+\delta)\Gamma^{\hat{t}}(\hat{t}+\delta)\bigg{]}+
      \frac{1}{\delta}\mathbb{E}\bigg{[}\int^{\hat{t}+\delta}_{{\hat{t}}}(\tilde{\cal{L}}{\varphi})(\hat{\gamma}_{{\hat{t}}},u^{\varepsilon}(l),\pi(u^{\varepsilon})(l))dl\bigg{]}\nonumber\\
                      &&+\frac{1}{\delta}\mathbb{E}\bigg{[}\int^{\hat{t}+\delta}_{{\hat{t}}}[(\tilde{\cal{L}}{\varphi})(X^{\hat{\gamma}_{\hat{t}},{u}^{\varepsilon},\pi({u}^{\varepsilon})}_l,u^{\varepsilon}(l),\pi(u^{\varepsilon})(l))-
                      (\tilde{\cal{L}}{\varphi})(\hat{\gamma}_{{\hat{t}}},u^{\varepsilon}(l),\pi(u^{\varepsilon})(l))]dl\bigg{]}\nonumber\\
                      &&+\frac{1}{\delta}\mathbb{E}\bigg{[}\int^{\hat{t}+\delta}_{{\hat{t}}}(\Gamma^{\hat{t}}(l)-1)(\tilde{\cal{L}}{\varphi})(X^{\hat{\gamma}_{\hat{t}},{u}^{\varepsilon},\pi({u}^{\varepsilon})}_l,
                                u^{\varepsilon}(l),\pi(u^{\varepsilon})(l))dl\bigg{]}\nonumber\\
    &:=&I+II+III+IV.
\end{eqnarray}
Since the coefficients in $\tilde{\cal{L}}$  satisfy  linear growth  condition,
 combining the regularity of $\varphi\in C^{1,2}_{p}({\Lambda}^{\hat{t}})$, there exist a integer
$\bar{p}\geq1$ and a constant $C>0$ independent of $(u,v)\in U\times V$ such that, for all $(t,\gamma_t,u,v)\in [0,T]\times\Lambda\times U\times V$,
\begin{eqnarray}\label{4.4444}|{\varphi}(\gamma_{t})|\vee  |
                      (\tilde{\cal{L}}{\varphi})(\gamma_{t},u,v)|
                      \leq  C(1+||\gamma_{t}||_0)^{\bar{p}}.
\end{eqnarray}
In view of Lemma \ref{lemmaexist}, we also have
\begin{eqnarray*}
                          \sup_{u^\varepsilon(\cdot)\in{\cal{U}}[\hat{t},\hat{t}+\delta]}\mathbb{E}[\sup_{{\hat{t}}\leq s\leq \hat{t}+\delta}|\Gamma^{\hat{t}}(s)-1|^2]\leq C\delta.
\end{eqnarray*}
Thus, by $\varphi\in {\cal{A}}^+(\hat{\gamma}_{\hat{t}},\underline{V})$,
\begin{eqnarray}\label{4.1611}
I&=& -\frac{1}{\delta}\mathbb{E}\bigg{[}\bigg{(}{\varphi}(X_{\hat{t}+\delta}^{\hat{\gamma}_{\hat{t}},{u}^{\varepsilon},\pi({u}^{\varepsilon})})-Y^{\hat{\gamma}_{\hat{t}},u^{\varepsilon},\pi(u^{\varepsilon})}(\hat{t}+\delta)\bigg{)}\Gamma^{\hat{t}}(\hat{t}+\delta)\bigg{]}\nonumber\\
&=&\frac{1}{\delta}\mathbb{E}\bigg{[}\bigg{(}\underline{V}(X^{\hat{\gamma}_{\hat{t}},
 {u}^{\varepsilon},\pi({u}^{\varepsilon})}_{{\hat{t}}+\delta})-{\varphi}(X_{\hat{t}+\delta}^{\hat{\gamma}_{\hat{t}},{u}^{\varepsilon},\pi({u}^{\varepsilon})})\bigg{)}\Gamma^{\hat{t}}(\hat{t}+\delta)\bigg{]}\leq 0;
\end{eqnarray}
and by (\ref{628jia}),
\begin{eqnarray}\label{4.16}
             II&\leq&\frac{1}{\delta}\bigg{[}\int^{\hat{t}+\delta}_{{\hat{t}}}\sup_{u\in U}(\tilde{\cal{L}}{\varphi})({\hat{\gamma}_{\hat{t}},{u},\pi(u)})dl\bigg{]}\nonumber\\
                      &=& \partial_t{\varphi}(\hat{\gamma}_{\hat{t}})
                           +\sup_{u\in U}{{\mathbf{H}}}(\hat{\gamma}_{\hat{t}},{\varphi}(\hat{\gamma}_{\hat{t}}),\partial_x{\varphi}(\hat{\gamma}_{\hat{t}}),
                           \partial_{xx}{\varphi}(\hat{\gamma}_{\hat{t}}),u,\pi(u))\leq -\theta.
\end{eqnarray}
  \par
Now we estimate higher order terms  $III$ and $IV$.  By (\ref{fbjia3}),
$$
\lim_{\delta\rightarrow0}\mathbb{E}d_\infty^{\bar{p}}(X^{\hat{\gamma}_{\hat{t}},{u}^{\varepsilon},\pi({u}^{\varepsilon})}_{\hat{t}+\delta},\hat{\gamma}_{\hat{t}})=0.
$$
Then by (\ref{4.4444}) 
   and the dominated convergence theorem, we have
$$
\lim_{\delta\rightarrow0} \mathbb{E}\sup_{\hat{t}\leq l\leq \hat{t}+\delta}|(\tilde{\cal{L}}{\varphi})(X^{{\hat{\gamma}_{\hat{t}},{u}^{\varepsilon},\pi({u}^{\varepsilon})}}_l,{u}^{\varepsilon}(l),\pi({u}^{\varepsilon})(l))-
                      (\tilde{\cal{L}}{\varphi})({\hat{\gamma}_{\hat{t}},{u}^{\varepsilon}(l),\pi({u}^{\varepsilon})}(l))|=0.
$$
Therefore,
\begin{eqnarray}\label{4.18}
\lim_{\delta\rightarrow0}|III| \leq\lim_{\delta\rightarrow0}\sup_{\hat{t}\leq l\leq \hat{t}+\delta}\mathbb{E}|(\tilde{\cal{L}}{\varphi})(X^{{\hat{\gamma}_{\hat{t}},{u}^{\varepsilon}(l),\pi({u}^{\varepsilon})}}_l,{u}^{\varepsilon}(l),\pi({u}^{\varepsilon})(l))-
                      (\tilde{\cal{L}}{\varphi})({\hat{\gamma}_{\hat{t}},{u}^{\varepsilon},\pi({u}^{\varepsilon})}(l))|=0;
\end{eqnarray}
and, for some finite constant $C>0$,
\begin{eqnarray}\label{4.19}
|IV|
&\leq&\frac{1}{\delta}\int^{\hat{t}+\delta}_{{\hat{t}}}\mathbb{E}|\Gamma^{\hat{t}}(l)-1||(\tilde{\cal{L}}{\varphi})(X^{\hat{\gamma}_{\hat{t}},{u}^{\varepsilon},\pi({u}^{\varepsilon})}_l,
                               u^{\varepsilon}(l),\pi(u^{\varepsilon})(l)|
                       dl\nonumber\\
                       &\leq&\frac{1}{\delta}\int^{\hat{t}+\delta}_{{\hat{t}}}(\mathbb{E}(\Gamma^{\hat{t}}(l)-1)^2)^{\frac{1}{2}}(\mathbb{E}((\tilde{\cal{L}}{\varphi})(X^{\hat{\gamma}_{\hat{t}},
                                     {u}^{\varepsilon},\pi({u}^{\varepsilon})}_l,
                               u^{\varepsilon}(l),\pi(u^{\varepsilon})(l))^2)^{\frac{1}{2}}dl
                    \nonumber\\
                       &\leq& C(1+||\hat{\gamma}_{\hat{t}}||_0)^{\bar{p}}\delta^\frac{1}{2}.
\end{eqnarray}
Substituting  (\ref{4.1611}), (\ref{4.16}) and (\ref{4.19}) into (\ref{4.15}), we have
\begin{eqnarray}\label{4.2000000}
-\varepsilon\leq -\theta+III+C(1+||\hat{\gamma}_{\hat{t}}||_0)^{\bar{p}}\delta^\frac{1}{2}.
\end{eqnarray}
Sending $\delta$ to $0$, and then $\varepsilon\rightarrow0$, by  (\ref{4.18}) we deduce that $\theta\leq 0$, which induces a contradiction. Therefore,
(\ref{628jia0}) holds true, and
by the arbitrariness of $\varphi\in {\cal{A}}^+(\hat{\gamma}_{\hat{t}},\underline{V})$, we show
 $\underline{V}$ is  a viscosity subsolution to (\ref{hjb1}).
 \par
 Now let us prove that $\underline{V}$ ia a viscosity supersolution of (\ref{hjb106}). We let  $\varphi\in {\cal{A}}^-(\hat{\gamma}_{\hat{t}},\underline{V})$
                  with
                   $(\hat{t},\hat{\gamma}_{\hat{t}})\in [0,T)\times \Lambda$.  
 For  $0< \delta\leq T-\hat{t}$, we have $\hat{t}< \hat{t}+\delta \leq T$, then by the DPP (Theorem \ref{theoremddp}), we obtain the following result:
 \begin{eqnarray}\label{4.92}
                           && 0=\underline{V}(\hat{\gamma}_{\hat{t}})-{{\varphi}} (\hat{\gamma}_{\hat{t}})
                           =\mathop{\essinf}\limits_{\beta(\cdot)\in {\cal{B}}_{[\hat{t},\hat{t}+\delta]}}\mathop{\esssup}\limits_{u(\cdot)\in{\mathcal
                              {U}}[\hat{t},\hat{t}+\delta]}G^{\hat{\gamma}_{\hat{t}},u,\beta(u)}_{\hat{t},\hat{t}+\delta}[\underline{V}(X^{\hat{\gamma}_{\hat{t}},u,\beta(u)}_{\hat{t}+\delta})]
                           -{{\varphi}} (\hat{\gamma}_{\hat{t}}).
\end{eqnarray}
Since $\mathop{\essinf}\limits_{v(\cdot)\in {\cal{V}}[\hat{t},\hat{t}+\delta]}G^{\hat{\gamma}_{\hat{t}},u,v}_{\hat{t},\hat{t}+\delta}[\underline{V}(X^{\hat{\gamma}_{\hat{t}},u,v}_{\hat{t}+\delta})]\leq G^{\hat{\gamma}_{\hat{t}},u,\beta(u)}_{\hat{t},\hat{t}+\delta}[\underline{V}(X^{\hat{\gamma}_{\hat{t}},u,\beta(u)}_{\hat{t}+\delta})]$ for all $u(\cdot)\in {\cal{U}}[\hat{t},\hat{t}+\delta]$ and $\beta(\cdot)\in {\cal{B}}_{[\hat{t},\hat{t}+\delta]}$,
 we get
 \begin{eqnarray*}
                           \mathop{\esssup}\limits_{u(\cdot)\in {\cal{U}}[\hat{t},\hat{t}+\delta]}\mathop{\essinf}\limits_{v(\cdot)\in {\cal{V}}[\hat{t},\hat{t}+\delta]}G^{\hat{\gamma}_{\hat{t}},u,v}_{\hat{t},\hat{t}+\delta}[\underline{V}(X^{\hat{\gamma}_{\hat{t}},u,v}_{\hat{t}+\delta})]
                           \leq
                             \mathop{\essinf}\limits_{\beta(\cdot)\in {\cal{B}}_{[\hat{t},\hat{t}+\delta]}}\mathop{\esssup}\limits_{u(\cdot)\in{\mathcal
                              {U}}[\hat{t},\hat{t}+\delta]}G^{\hat{\gamma}_{\hat{t}},u,\beta(u)}_{\hat{t},\hat{t}+\delta}[\underline{V}(X^{\hat{\gamma}_{\hat{t}},u,\beta(u)}_{\hat{t}+\delta})].
 \end{eqnarray*}
Then, for any $u\in U$, $\varepsilon>0$ and $0<\delta\leq T-\hat{t}$,  we can  find a control  ${v}^{\varepsilon}(\cdot)
\equiv v^{{\varepsilon},\delta}(\cdot)
\in {\cal{V}}[\hat{t},\hat{t}+\delta]$ such
   that 
   the following result holds:
\begin{eqnarray}\label{4.100}
     G^{\hat{\gamma}_{\hat{t}},{u},{v}^{\varepsilon}}_{{\hat{t}},\hat{t}+\delta}[\underline{V}(X^{\hat{\gamma}_{\hat{t}},{u},{v}^{\varepsilon}}_{{\hat{t}}+\delta})]-{{\varphi}} (\hat{\gamma}_{\hat{t}})
     \leq  {\varepsilon}\delta.
\end{eqnarray}
Now following similar arguments as above we can show
$$
\partial_t{\varphi}(\hat{\gamma}_{\hat{t}})+{\mathbf{H}}^-(\hat{\gamma}_{\hat{t}},{\varphi}(\hat{\gamma}_{\hat{t}}),\partial_x{\varphi}(\hat{\gamma}_{\hat{t}}),
                           \partial_{xx}{\varphi}(\hat{\gamma}_{\hat{t}}))
\leq0.
$$
By the arbitrariness of $\varphi\in {\cal{A}}^-(\hat{\gamma}_{\hat{t}},\underline{V})$, we show
 $\underline{V}$ is  a viscosity supersolution to (\ref{hjb1}).
                 This step  completes the proof.\ \ $\Box$
                   \par
 To complete the previous proof, it remains to state and prove the following lemma.
 \\
 {\bf Lemma \ \ B.1.} \ \
                 Assume that Hypothesis \ref{hypstate6}  holds. Then there exists a measurable function $\pi:U\rightarrow V$ such that (\ref{628jia}) holds true.

\par
   {\bf  Proof}. \ \  Define $\Pi:U\rightarrow 2^V$ by
   $$
                       \Pi(u)=\{v\in V|\partial_t{\varphi}(\hat{\gamma}_{\hat{t}})
                           +\mathbf{H}(\hat{\gamma}_{\hat{t}},{\varphi}(\hat{\gamma}_{\hat{t}}),\partial_x{\varphi}(\hat{\gamma}_{\hat{t}}),\partial_{xx}{\varphi}(\hat{\gamma}_{\hat{t}}),u,v)\leq-\theta\}, \ \  u\in U.
   $$
   It is clear that $\Pi(u)$ is a nonempty and closed subset. Now we prove $\Pi$ is measurable. For every compact set $O\in V$, define
   $$
   \Pi^-(O)=\{u\in U|\Pi(u)\cap O\neq\phi\}.
   $$
   Assume $\{u_i\}_{i\geq1}\in \Pi^-(O)$ and $u_i\rightarrow u$ in $U$, then there exist $\{v_i\}_{i\geq1}\in V$ such that
\begin{eqnarray}\label{pi}
   v_i\in  O,\ \mbox{and} \ \ \partial_t{\varphi}(\hat{\gamma}_{\hat{t}})
                           +\mathbf{H}(\hat{\gamma}_{\hat{t}},{\varphi}(\hat{\gamma}_{\hat{t}}),\partial_x{\varphi}(\hat{\gamma}_{\hat{t}}),\partial_{xx}{\varphi}(\hat{\gamma}_{\hat{t}}),u_i,v_i)\leq-\theta.
\end{eqnarray}
   Since $O$ is a compact subset, then a $v\in O$ and a  subsequence of $\{v_i\}_{i\geq1}$ still denoted by itself  exist such that $v_i\rightarrow v$ in $V$.  Letting $i\rightarrow\infty$ in (\ref{pi}), we get that
   $$
    \partial_t{\varphi}(\hat{\gamma}_{\hat{t}})
                           +\mathbf{H}(\hat{\gamma}_{\hat{t}},{\varphi}(\hat{\gamma}_{\hat{t}}),\partial_x{\varphi}(\hat{\gamma}_{\hat{t}}),\partial_{xx}{\varphi}(\hat{\gamma}_{\hat{t}}),u,v)\leq-\theta.
   $$
   Thus, $u\in \Pi^-(O)$. Therefore, $\Pi^-(O)$ is closed. Then, by Proposition 4.2.9 in \cite{den}, $\Pi$ is measurable. From Theorem 4.3.1 in \cite{den}, it follows that $\Pi$ admits a  measurable selection $\pi:U\rightarrow V$ such that (\ref{628jia}) holds true.  The proof is now complete.
 \ \ $\Box$

\par

\end{document}